\documentclass[12pt,a4paper]{article}
\usepackage{latexsym,epsfig,fullpage,amsfonts,amssymb,amsmath,amscd,epic,graphics,
rotating,theorem,yhmath,pxfonts,comment,tikz}
\usepackage[all]{xy}
\theoremstyle{plain}
\newtheorem{lemma}{Lemma}[section]
\newtheorem{proposition}[lemma]{Proposition}
\newtheorem{remark}[lemma]{Remark}
\newtheorem{example}[lemma]{Example}
\newtheorem{theorem}[lemma]{Theorem}
\newtheorem{definition}[lemma]{Definition}
\newtheorem{corollary}[lemma]{Corollary}
{\theorembodyfont{\rmfamily}

\newcommand{\oa}{\overline{a}}

\tikzstyle{regular_node}=[circle, draw, fill=black!50,
                        inner sep=0pt, minimum width=4pt]
\tikzstyle{small_node}=[circle, draw, fill=black!50,
                        inner sep=0pt, minimum width=2pt]

\begin{document}
\newcommand{\pperp}{\hbox{$\perp\hskip-6pt\perp$}}
\newcommand{\ssim}{\hbox{$\hskip-2pt\sim$}}
\newcommand{\aleq}{{\ \stackrel{3}{\le}\ }}
\newcommand{\ageq}{{\ \stackrel{3}{\ge}\ }}
\newcommand{\aeq}{{\ \stackrel{3}{=}\ }}
\newcommand{\bleq}{{\ \stackrel{n}{\le}\ }}
\newcommand{\bgeq}{{\ \stackrel{n}{\ge}\ }}
\newcommand{\beq}{{\ \stackrel{n}{=}\ }}
\newcommand{\cleq}{{\ \stackrel{2}{\le}\ }}
\newcommand{\cgeq}{{\ \stackrel{2}{\ge}\ }}
\newcommand{\ceq}{{\ \stackrel{2}{=}\ }}
\newcommand{\A}{{\mathbb A}}\newcommand{\bL}{{\boldsymbol L}}
\newcommand{\K}{{\mathbb K}}
\newcommand{\Z}{{\mathbb Z}}\newcommand{\F}{{\mathbf F}}
\newcommand{\R}{{\mathbb R}}
\newcommand{\C}{{\mathbb C}}
\newcommand{\Q}{{\mathbb Q}}
\newcommand{\PP}{{\mathbb P}}
\newcommand{\mnote}{\marginpar}\newcommand{\red}{{\operatorname{red}}}
\newcommand{\Leaf}{{\operatorname{Leaf}}}\newcommand{\Sym}{{\operatorname{Sym}}}
\newcommand{\oeps}{{\overline\eps}}\newcommand{\Div}{{\operatorname{Div}}}
\newcommand{\oDel}{{\widetilde\Del}}\newcommand{\lab}{{\operatorname{lab}}}
\newcommand{\real}{{\operatorname{Re}}}\newcommand{\Aut}{{\operatorname{Aut}}}
\newcommand{\conv}{{\operatorname{conv}}}\newcommand{\BG}{{\operatorname{BG}}}
\newcommand{\Span}{{\operatorname{Span}}}\newcommand{\GS}{{\operatorname{GS}}}
\newcommand{\Ker}{{\operatorname{Ker}}}\newcommand{\Ss}{{\operatorname{SS}}}
\newcommand{\Ann}{{\operatorname{Ann}}}\newcommand{\mt}{{\operatorname{wt}}}
\newcommand{\Fix}{{\operatorname{Fix}}}\newcommand{\Ima}{{\operatorname{Im}}}
\newcommand{\sign}{{\operatorname{sign}}}\newcommand{\Def}{{\operatorname{Def}}}
\newcommand{\Tors}{{\operatorname{Tors}}}
\newcommand{\Card}{{\operatorname{Card}}}
\newcommand{\alg}{{\operatorname{alg}}}\newcommand{\ord}{{\operatorname{ord}}}
\newcommand{\oi}{{\overline i}}
\newcommand{\oj}{{\overline j}}
\newcommand{\bb}{\boldsymbol{b}}
\newcommand{\os}{{\overline s}}
\newcommand{\ba}{\boldsymbol{a}}
\newcommand{\be}{\boldsymbol{e}}
\newcommand{\ow}{{\overline w}}
\newcommand{\bc}{\boldsymbol{c}}\newcommand{\om}{\overline{m}}
\newcommand{\oz}{{\overline z}}\newcommand{\on}{{\overline n}}
\newcommand{\eps}{{\varepsilon}}
\newcommand{\proofend}{\hfill$\Box$\bigskip}
\newcommand{\Int}{{\operatorname{Int}}}\newcommand{\grad}{{\operatorname{grad}}}
\newcommand{\pr}{{\operatorname{pr}}}
\newcommand{\Hom}{{\operatorname{Hom}}}
\newcommand{\Ev}{{\operatorname{Ev}}}
\newcommand{\im}{{\operatorname{Im}}}\newcommand{\br}{{\operatorname{br}}}
\newcommand{\sk}{{\operatorname{sk}}}\newcommand{\Fl}{{\operatorname{Fl}}}
\newcommand{\const}{{\operatorname{const}}}
\newcommand{\Sing}{{\operatorname{Sing}}\hskip0.06cm}
\newcommand{\conj}{{\operatorname{Conj}}}
\newcommand{\Cl}{{\operatorname{Cl}}}
\newcommand{\Crit}{{\operatorname{Crit}}}
\newcommand{\Ch}{{\operatorname{Ch}}}
\newcommand{\discr}{{\operatorname{discr}}}
\newcommand{\Tor}{{\operatorname{Tor}}}
\newcommand{\Conj}{{\operatorname{Conj}}}
\newcommand{\vol}{{\operatorname{vol}}}
\newcommand{\defect}{{\operatorname{def}}}
\newcommand{\codim}{{\operatorname{codim}}}
\newcommand{\tmu}{{\C\mu}}
\newcommand{\bv}{\boldsymbol{v}}
\newcommand{\ox}{{\overline{x}}}
\newcommand{\bw}{{\boldsymbol w}}\newcommand{\bn}{{\boldsymbol n}}
\newcommand{\bx}{{\boldsymbol x}}
\newcommand{\bd}{{\boldsymbol d}}
\newcommand{\bz}{{\boldsymbol z}}\newcommand{\bp}{{\boldsymbol p}}
\newcommand{\tet}{{\theta}}
\newcommand{\Del}{{\Delta}}
\newcommand{\bet}{{\beta}}
\newcommand{\kap}{{\kappa}}
\newcommand{\del}{{\delta}}
\newcommand{\sig}{{\sigma}}
\newcommand{\alp}{{\alpha}}
\newcommand{\Sig}{{\Sigma}}
\newcommand{\Gam}{{\Gamma}}
\newcommand{\gam}{{\gamma}}\newcommand{\idim}{{\operatorname{idim}}}
\newcommand{\Lam}{{\Lambda}}
\newcommand{\lam}{{\lambda}}
\newcommand{\SC}{{SC}}
\newcommand{\MC}{{MC}}
\newcommand{\nek}{{,...,}}
\newcommand{\cim}{{c_{\mbox{\rm im}}}}
\newcommand{\clM}{\tilde{M}}
\newcommand{\clV}{\bar{V}}

\title{On refined count of rational tropical curves}
\author{Eugenii Shustin\thanks{School of Mathematical Sciences, Tel Aviv University,
Ramat Aviv, 69978 Tel Aviv, Israel. E-mail: shustin@tauex.tau.ac.il}}
\date{}
\maketitle

\centerline{\it\small To my good friend and a remarkable mathematician Gert-Martin Greuel}
\centerline{\it\small on occasion of his 75th
birthday}

\begin{abstract}
We address the problem of existence of refined (i.e., depending on a formal parameter) tropical enumerative invariants, and we present two new examples of a refined count of
rational marked tropical curves.
One of the new invariants counts plane rational tropical curves with an unmarked vertex of arbitrary valency. It was motivated by the tropical enumeration of plane cuspidal tropical curves given by Y. Ganor and the author, which naturally led to consideration of
plane tropical curves with an unmarked four-valent vertex. Another refined invariant counts rational tropical curves of a given degree in the Euclidean space of arbitrary dimension matching specific constraints, which make the spacial refined invariant similar to known planar invariants.
\end{abstract}

{\it MSC-2010}: 14N10, 14T05

\section*{Introduction}
Refined enumerative invariants, i.e., invariants that count objects with weights depending on a parameter, appeared in
\cite{GoS}. Later F. Block and L. G\"ottsche \cite{BG} introduced refined invariants
in the framework of tropical geometry.
This subject has been developed further in \cite{BS,BG1,GS,IM,Ma,MaRu1,SS}.
The main challenging problem in this direction is to find a reasonable
algebraic-geometric or symplectic
counterpart of the refined invariants. Specific values of some of these refined invariants
coincide with closed or open Gromov-Witten invariants of toric surfaces, or with descendant invariants of toric surfaces. In a certain particular case the relation
between the refined tropical and refined algebraic-geometric count
was discovered by G. Mikhalkin \cite{Mi4}. A related important question is: Which numerical tropical enumerative invariants admit a refinement?

In this note, we present
two new examples of a refined count of
rational marked tropical curves.

One of the new invariants counts plane rational tropical curves with an unmarked vertex of arbitrary valency. It was motivated by the correspondence between plane unicuspidal algebraic curves and
plane cuspidal tropical curves given by Y. Ganor and the author \cite{GSh}, where we encountered
plane tropical curves with an unmarked four-valent vertex. The weight of such a cuspidal tropical curve appeared to be convenient for a refinement, while the refined invariant specializes to the number of plane rational unicuspidal curves of a given degree that match an appropriate point constraint, as the parameter gets value $1$.

Another refined invariant counts rational tropical curves of a given degree in the Euclidean space of arbitrary dimension.
The known planar refined invariants \cite{BS,BG,GS} come from numerical tropical enumerative invariants, which count plane tropical curves of a given degree and genus
so that the weight assigned to each tropical curve is
the product of weights of vertices of this tropical curve.
Then the weight of each vertex is replaced with its refined analogue,
and the final task is to show that this formally obtained expression
does not depend on the choice of constraints.
For tropical curves in $\R^n$, $n\ge3$, the first tropical
enumerative invariant was
found by T. Nishinou and B. Siebert \cite{NS}. It then was largely extended by T. Mandel and H. Ruddat
\cite{MaRu} by allowing $\psi$-constraints, complicated boundary conditions, and curves of positive genera.
The weights of tropical curves in these enumerative problems involve factors depending
on the entire curves, not only on their vertices, edges or other small
fragments. This, in fact, makes difficult to find an appropriate
refinement of the numerical invariant. Our main idea is to specify a
particular situation when the weight of each tropical curve
splits into relatively elementary factors, which then are replaced by their refined analogues.
Independently, a somewhat similar idea was elaborated in \cite{MaRu1}, where the authors use collections of constraints different from ours and also allow multivalent marked points.

\smallskip
{\bf Acknowledgements}. The author has been supported by the Israel Science Foundation grants no. 176/15 and 501/18, and by The Bauer-Neuman Chair in Real and Complex Geometry. The main part of this research was performed during the author's stay at the Max-Planck Institut f\"ur Mathematik in August-September 2018.
The author thanks the MPIM for support and excellent working conditions. Special thanks are due to T. Mandel and H. Ruddat, who attracted my attention to their works \cite{MaRu,MaRu1} and explained me important details of their approach both to numerical and refined tropical enumerative invariants.
I am also grateful to T. Blomme, M. Polyak, and the unknown referee for valuable remarks which helped me to correct mistakes and improve the presentation.

\section{Marked rational tropical curves}
We shortly recall some basic definitions concerning rational tropical
curves adapted to our setting and define
the class of tropical curves under consideration (for details,
see \cite{GKM,GM,Mi,Mi3}).

\smallskip

{\bf(1)} A {\it $n$-marked rational tropical curve in $\R^m$},
$m\ge2$, is a triple $(\Gamma,h,\bp)$, where
\begin{itemize}\item $\Gamma$ is a finite metric graph, which
is either isometric to $\R$, or is a finite
connected metric tree without vertices of valency $\le2$, whose set $\Gamma^0$ of vertices in nonempty, the
set of edges $\Gamma^1$ contains a
subset $\Gamma^1_\infty\ne\emptyset$ consisting of edges isometric to $[0,\infty)$ (called {\it ends}), while
$\Gamma^1\setminus\Gamma^1_\infty$ consists of edges isometric to compact segments in $\R$ (called {\it finite edges});
\item $h:\Gamma\to\R^m$ is a proper continuous map such that
$h$ is nonconstant, affine-integral on each edge of $\Gamma$ in the length coordinate
and, at each vertex $V$ of $\Gamma$, the {\it balancing condition} holds
\begin{equation}\sum_{E\in\Gamma^1,\ V\in E}
\ba_V(E)=0\ ,\label{se3}\end{equation} where $\ba_V(E)$
is the image under the differential $D(h\big|_E)$
of the unit tangent vector to $E$ emanating from
its endpoint $V$ (called the {\it directing vector} of $E$ centered at $V$);
\item $\bp=(p_1,...,p_n)$ is a sequence of $n$ points of $\Gamma$.
\end{itemize}

The multiset of vectors $\deg(\Gamma,h):=\big\{\ba_V(E)\ :\ E\in\Gamma^1_\infty\big\}\subset\Z^m\setminus\{0\}$ is called the {\it degree} of $(\Gamma,h,\bp)$.
Clearly the vectors of $\deg(\Gamma,h)$ sum up to zero (we call such a multiset {\it balanced}).
The degree $\Del$ is called {\it nondegenerate} if $\dim\Span\{\ba\in\Del\}=m$, and is called
{\it primitive} if all vectors $\ba\in\Del$ are primitive integral vectors.

\smallskip

{\bf(2)} Two plane $n$-marked rational tropical curves $(\Gamma,h,\bp)$ and $(\Gamma',h',\bp')$ are called isomorphic, if there exists
an isometry $\varphi:\Gamma\to\Gamma'$ such that $h=h'\circ\varphi$ and $\varphi(p_i)=p'_i$ for all
$p_i\in\bp$, $p'_i\in\bp'$, $i=1,...,n$. Clearly, the isomorphism preserves the {\it combinatorial type} of a tropical curve $(\Gamma,h,\bp)$, i.e., the combinatorial type of the pair $(\Gamma,\bp)$ and the list of all
directing vectors $\ba_V(E)$ for all pairs of incident vertices and edges of $\Gamma$. For a given degree $\Del$, there are only finitely many
combinatorial types of pairs $(\Gamma,\bp)$ (see \cite[Proposition 2.1]{NS}).
Given a finite balanced multiset
$\Del\subset\Z^m\setminus\{0\}$, we consider the moduli space ${\mathcal M}_{0,n}(\R^m,\Del)$ parameterizing
isomorphism classes $[(\Gamma,h,\bp)]$ of $n$-marked rational tropical curves of degree $\Del$ in $\R^m$.

We will also use labeled tropical curves. For, we fix a linear order on $\Del$ and denote the obtained sequence by $\Del^\lab$. A {\it labeled $n$-marked plane rational tropical curve of degree}
$\Del^\lab$ is a triple $(\Gam^\lab,\bp,h)$, where $(\Gam,\bp,h)$ is an $n$-marked plane
rational tropical curve of degree $\Del$, and $\Gam^\lab$ is the graph $\Gam$ with a linear order on the set of its ends $\Gam^1_\infty$ such that $h$ equips the $i$-th end of $\Gam^\lab$ with the
$i$-th directing vector in $\Del^\lab$ for all $i$
(cf. \cite[Definition 4.1]{GKM}). Note that a labeled tropical curve has no nontrivial isomorphisms, and hence the corresponding moduli space ${\mathcal M}^\lab_{0,n}(\R^m,\Delta)$ parameterizes just the $n$-marked labeled rational tropical curves of degree $\Del^\lab$.

\begin{lemma}[\cite{GKM}, Proposition 7.4, and \cite{BS}, Lemma 1.1]\label{sl1}
The space ${\mathcal M}^\lab_{0,n}(\R^m,\Del)$ can be identified with a tropical fan of dimension
$|\Del|+m+n-3$ in $\R^N$ for a sufficiently large $N$, whose open cells parameterize isomorphism classes of marked tropical curves  of the same combinatorial type. The map forgetting labels
$$\pi_{0,n}:{\mathcal M}^{\lab}_{0,n}(\R^m,\Del^\lab)\to{\mathcal M}_{0,n}(\R^m,\Delta)$$
is surjective and finite, and, for any element $[(\Gamma,\bp,h)]\in
{\mathcal M}_{0,n}(\R^m,\Delta)$, we have
\begin{equation}\left|(\pi_{0,n})^{-1}(\Gamma,h,\bp)\right|=\frac{|G|!}{|\Aut(\Gamma,h,\bp)|}\ ,
\label{ce5}\end{equation}
where $G$ is the group of permutations $\varphi:\Del^\lab\to\Del^\lab$ such that $\varphi(\bb)=\bb$ for each $\bb\in\Del^\lab$.
\end{lemma}

\section{The refined cuspidal tropical invariant of the plane}\label{csec1}
In this section we consider only plane tropical curves.

\subsection{Moduli spaces of marked plane tropical curves with multi-valent vertices}
Denote by $\Z^\infty_+$ the set of sequences of nonnegative integers $(k_i)_{i\ge0}$ such that $\sum_ik_i<\infty$, and by $\Z^{\infty,*}_+\subset\Z^\infty_+$
the set of sequences with the vanishing initial member.
Let $\om=(m_i)_{i\ge0}\in\Z^{\infty,*}_+$, $\on=(n_i)_{i\ge0}\in\Z^\infty_+$.
We say that a marked plane tropical curve $(\Gamma,h,\bp)$ is of V-type $(\om,\on)$,  if exactly $m_i$ vertices in $\Gamma^0\setminus\bp$ have valency $i+2$ for all $i\ge1$,
exactly $n_0$ points of $\bp$ lie in $\Gamma\setminus\Gamma^0$,
exactly $n_i$ vertices in $\Gamma^0\cap\bp$ have valency $i+2$ for all $i\ge1$. It is easy to see that
\begin{equation}|\Gamma^0|=\sum_{i\ge1}(m_i+n_i),\quad|\Gamma^1|=\sum_{i\ge1}(i+1)
(m_i+n_i)+1,\quad|\Gamma^1_\infty|=
\sum_{i\ge1}i(m_i+n_i)+2\ .\label{ce1}\end{equation}
Assuming that a finite balanced multiset
$\Delta\subset\Z^2\setminus\{0\}$ and sequences of nonnegative integers
$\om\in\Z^{\infty,*}_+$, $\on\in\Z^\infty_+$, satisfy
\begin{equation}|\Delta|=\sum_{i\ge1}i(m_i+n_i)+2\label{ce3}\end{equation} (cf. (\ref{ce1})), we consider the moduli space ${\mathcal M}_{0,\on,\om}(\R^2,\Delta)$ parameterizing
isomorphism classes $[(\Gamma,h,\bp)]$ of plane $n=\sum_{i\ge0}n_i$-marked rational tropical curves of V-type $(\om,\on)$ and degree $\Delta$.

\begin{lemma}\label{cl1}
The space ${\mathcal M}_{0,\on,\om}(\R^2,\Delta)$ can be identified with a finite union of open convex polyhedral cones of
pure dimension $\sum_{i\ge0}(m_i+n_i)+1$.
\end{lemma}

{\bf Proof.}
Given a combinatorial type of the pair $(\Gamma,\bp)$ and the distribution of the directing vectors $\ba_V(E)\in\Z^2\setminus\{0\}$
for all edges $E\in\Gamma^1$, the lengths of the finite edges, the distances from marked points in $\Gamma\setminus\Gamma^0$ to chosen vertices of the corresponding edges, and
the freely chosen image $h(V)$ of a fixed vertex $V\in\Gamma^0$ give $N=\sum_{i\ge0}(m_i+n_i)+1$  independent coordinates in the positive orthant of
$\R^N$, from which one should get rid suitable diagonals in case when more than one marked points occur in the interior of the same edge of $\Gamma$.
\proofend

By $\widehat{\mathcal M}_{0,\on,\om}(\R^2,\Delta)$ we denote the polyhedral fan obtained by extending ${\mathcal M}_{0,\on,\om}(\R^2,\Delta)$
with classes corresponding either to contraction of some finite edges (i.e., vanishing of their lengths), or to arrival of marked points from edges to
vertices of $\Gamma$, or to collision of marked points.

Assume that $2n=\sum_{i\ge0}(m_i+n_i)+1$, or, equivalently,
\begin{equation}n=\sum_{i\ge1}m_i+1\ .\label{ce2}\end{equation}
Then the evaluation map
$$\Ev:\widehat{\mathcal M}_{0,\on,\om}(\R^2,\Delta)\to\R^{2n},\quad \Ev[(\Gamma,h,\bp)]=h(\bp)\in\R^{2n}\ ,$$
relates spaces of the same dimension $\sum_{i\ge0}(m_i+n_i)+1=2n$.

\begin{definition}\label{cd1}
Let a balanced, nondegenerate multiset $\Delta\subset\Z^2\setminus\{0\}$ and sequences $\om\in\Z^{\infty,*}_+$,
$\on\in\Z^\infty_+$ satisfy (\ref{ce3}) and (\ref{ce2}).

(1) We say that a class $[(\Gamma,h,\bp)]\in{\mathcal M}_{0,\on,\om}(\R^2,\Delta)$ is regular, if
each connected component $K$ of $\Gamma\setminus\bp$ is unbounded, and its closure $\overline K\subset\Gamma$
possesses a unique orientation of its edges (called regular orientation) such that
\begin{itemize}\item all marked points in $\overline K$ are sources, all ends of $\overline K$ are oriented towards infinity;
\item for each vertex $V\in K\cap\Gamma^0$ exactly two of its incident edges are incoming, and, moreover, the $h$-images of these edges are
not collinear.
\end{itemize}

(2) A cell of ${\mathcal M}_{0,\on,\om}(\R^2,\Delta)$ is called enumeratively essential, if $\Ev$ injectively takes it to $\R^{2n}$. Denote by
${\mathcal M}^{\;e}_{0,\on,\om}(\R^2,\Delta)$ the union of the enumeratively essential cells of ${\mathcal M}_{0,\on,\om}(\R^2,\Delta)$, by
$\widehat{\mathcal M}^{\;e}_{0,\on,\om}(\R^2,\Delta)$ the closure of ${\mathcal M}^{\;e}_{0,\on,\om}(\R^2,\Delta)$ in
$\widehat{\mathcal M}_{0,\on,\om}(\R^2,\Delta)$, and by $\Ev^e$ the restriction of $\Ev$ to $\widehat{\mathcal M}^{\;e}_{0,\on,\om}(\R^2,\Delta)$.
\end{definition}

\begin{lemma}\label{cl3}
Let $\Delta\subset\Z^2\setminus\{0\}$ be a balanced nondegenerate multiset, $\om\in\Z^{\infty,*}_+$,
$\on\in\Z^\infty_+$. Suppose that
(\ref{ce3}) and (\ref{ce2}) hold. Then ${\mathcal M}^{\;e}_{0,\on,\om}(\R^2,\Delta)\ne\emptyset$ and
each cell of
${\mathcal M}^{\;e}_{0,\on,\om}(\R^2,\Delta)$ consists of regular classes.
\end{lemma}

{\bf Proof.} Suppose that $[(\Gamma,h,\bp)]\in{\mathcal M}_{0,\on,\om}(\R^2,\Delta)$ is a regular class. Then it belongs to ${\mathcal M}^{\;e}_{0,\on,\om}(\R^2,\Delta)$. Indeed,
whenever we fix the position of $h(\bp)$, the image of $h:\Gamma\to\R^2$ is fixed as well (recall that the combinatorial type of $(\Gamma,\bp)$
and the differentials of $h$ on the edges of $\Gamma$ are a priori fixed), and then we recover the lengths of
compact edges of $\Gamma$. On the other hand, it immediately follows from the regularity that any small variation of $h(\bp)$ induces a (unique) small variation
of $(\Gamma,h,\bp)$ in its combinatorial class.

For the proof of the existence of a regular class, we make the following elementary observation
(left to the reader as an exercise):

(O) Let $|\Delta|>3$ and let $\Delta$ by cyclically ordered by rotation in the positive direction. Then, for any $\ba_i\in\Delta$, which is not simultaneously collinear to $\ba_{i-1}$ and $\ba_{i+1}$, and any $1\le j\le |\Delta|-2$, there exist a sequence $\ba_k,...,\ba_{k+j}\in\Delta$,
including $\ba_i$, such that
$\dim\Span\{\ba_k,...,\ba_{k+j}\}=2$, $\sum_{s=k}^{k+j}\ba_s\ne0$, and the multiset
$\Delta'=(\Delta\setminus\{\ba_k,...,\ba_{k+j}\})\cup\{\ba'\}$, where $\ba'=\ba_k+...+\ba_{k+j}$, is balanced and nondegenerate.

\smallskip
Then we proceed in the same way as
the existence statement in \cite[Lemma 1.4]{BS}. We remind here this argument referring to \cite{BS} for the details.
First, we construct the (convex lattice) Newton polygon $P(\Delta)$, whose boundary can be represented as the union of
cyclically ordered integral segments $[v_k,v_{k+1}]$, $k=1,...,|\Delta|$,
$v_{|\Delta|+1}=v_1$, obtained by
rotating the ordered as above vectors $\ba_k\in\Delta$, $k=1,...,|\Delta|$, by $\frac{\pi}{2}$ clockwise (we call $\ba_k$ and $[v_k,v_{k+1}]$ dual to each other).
The set ${\mathcal V}=\{v_1,...,v_{|\Delta|}\}$ includes all the vertices of $P(\Delta)$.

We proceed by induction on $n$. If $n=1$, the curve $\Gamma$ is a fan with the center at the unique (marked) vertex, and $(\Gamma,h,\bp)$ is regular.
Suppose that $n>1$. Then (cf. (\ref{ce2})) there are $n_i>0$ and $m_j>0$. If $j=|\Delta|-2$ and respectively $n=n_0=2$,
then $\gamma$ again has a unique (unmarked) vertex, and we pick two marked points on two ends with non-collinear directing vectors, obtaining
a regular curve $(\Gamma,h,\bp)$. Suppose that $j\le|\Delta|-3$. If $i=0$, we choose $\ba_k$, which is not parallel both to $\ba_{k-1}$ and $\ba_{k+1}$,
then find a sequence $\ba_s,...,\ba_{s+j-1}$ as in observation (O), and, finally draw the chord in $P(\Delta)$ joining the points
$v_s$ and $v_{s+j}$. It follows that the interior of this chord is disjoint from $\partial P(\Delta)$. The chord cuts $P(\Delta)$ into a polygon
containing $j+2$ points of ${\mathcal V}$ and the remaining polygon $P(\Delta')$, where $\Delta'=(\Delta\setminus\{\ba_s,...,\ba_{s+j-1}\})\cup\{\ba'\}$,
$\ba'=\ba_s+...+\ba_{s+j-1}$. The former polygon (with the corresponding part of ${\mathcal V}$) is dual to a tropical curve with the unique
(unmarked) vertex of valency $j+2$, a marked point on the end directed by the vector $\ba_k$, and the end directed by the vector $-\ba'$, to which we attach
the remaining part of the constructed curve existing due to the induction assumption applied to $\Delta'$ and $\on'$, $\om'$, obtained
by reducing $n_i$ and $m_j$ by one. The regularity of the constructed tropical curve is evident. Suppose that $n_0=0$ and $i>0$. Then $i+j\le|\Delta|-3$.
We, first, choose a sequence $\ba_k,...,\ba_{k+i}$ as in observation (O), join the points $v_k,v_{k+i}\in{\mathcal V}$ by a chord, whose interior must be disjoint from
$\partial P(\Delta)$. It cuts off $P(\Delta)$ an polygon $P_1$ that will be dual to a marked point of valency $i+2$ incident to $i+1$ ends directed by
$\ba_k,...,\ba_{k+i}$ and to one more edge dual to the chord. Set $\Delta'=(\Delta\setminus\{\ba_k,...,\ba_{k_i}\})\cup\{\ba'\}$,
$\ba'=\ba_k+...,\ba_{k+i}$. Since the chord is not collinear with the neighboring sides of $P(\Delta)$, we apply observation (O) to $\Delta'$ and obtain a sequence
of $j+1$ vectors of $\Delta'$ (including $\ba'$), whose dual segments form a connected part of $\partial P(\Delta')$, and the extreme points of this part
are joined by a chord which intersects $\partial P(\Delta')$ only in its endpoints. Thus, we cut off $P(\Delta')$ a polygon $P_2$ which will be dual to an unmarked
vertex of valency $j+2$ incident to $j$ ends, an edge dual to the first constructed chord, and one more edge dual to the second chord. So, we attach the two constructed fragment by gluing along the edges dual to the first chord, and, finally, apply the induction assumption to $\Delta''$ that is formed by the vectors
$\ba_s$ dual to the remaining segments $[v_r,v_{r+1}]$ and by the vector $\ba''$ equal to the sum of all removed vectors $\ba_s$ (and dual to the second chord),
while $\on$ and $\om$ turn into $\on',\om'$ by reducing $1$ from $n_i$ and $m_j$.

\smallskip
Suppose now that $[(\Gamma,h,\bp)]\in{\mathcal M}^{\;e}_{0,\on,\om}(\R^2,\Delta)$. This means that $\bx=h(\bp)$ is a general position in $\R^2$.
Using induction on $|\Gamma^0|$, we show that $(\Gamma,h,\bp)$ is regular.
If $|\Gamma^0|\le1$, the claim is evident.
Assume that $|\Gamma^0|>1$. Since one can consider all components of $\Gamma\setminus\bp$ separately and independently, we are left with the case
when all marked points belong to ends of $\Gamma$, when no two points lie on the same end or on collinear ends. The relation $n=|\Gamma^0|+1$
(cf. (\ref{ce2})) yields that
there are two ends with marked points incident to the same vertex $V\in\Gamma^0$. Note that no any other end with a marked point is incident to $V$
due to the general position of $\bx$. So, we orient the segments on the chosen above two ends of $\Gamma$, which join the marked points with $V$, towards $V$,
while all other edges of $\Gamma$ incident to $V$ are oriented outwards. Thus, we reduce the considered case to the study of the connected components
of $\Gamma\setminus\{V\}$, and hence derive the required regularity by the induction assumption.
\proofend

\begin{remark}\label{cr2}
As shown in \cite[Proof of Lemma 1.4]{BS}, the subdivision of $P(\Delta)$ constructed in the proof of Lemma \ref{cl3} can be further refined
by extra chords between the points of ${\mathcal V}$ so that the final subdivision will consist of $|\Delta|-2$
nondegenerate triangles with vertices in ${\mathcal V}$.
\end{remark}

Let
\begin{equation}Y^{2n-1}=\Ev\big(\widehat{\mathcal M}^{\;e}_{0;\om}(\R^2,\Delta)\setminus{\mathcal M}^{\;e}_{0;\om}(\R^2,\Delta)\big)
\ .\label{ce11}\end{equation} This set possesses a structure of a finite polyhedral complex of dimension $\le2n-1$ in $\R^{2n}$, induced by that of
$\widehat{\mathcal M}^{\;e}_{0;\om}(\R^2,\Delta)\setminus{\mathcal M}^{\;e}_{0;\om}(\R^2,\Delta)$. Indeed, the images of (finitely many) cells of
$\widehat{\mathcal M}^{\;e}_{0;\om}(\R^2,\Delta)\setminus{\mathcal M}^{\;e}_{0;\om}(\R^2,\Delta)$ are polyhedra of dimension $\le2n-1$. Each of these polyhedra is an intersection of finitely many half-spaces. Let $H_1,...,H_s\subset\R^{2n}$ be the supporting hyperplanes of all the half-spaces that appear here,
and let $H^+_i,H^-_i$ be the half-spaces supported by $H_i$, $i=1,...,s$. For each point $\bx\in\R^{2n}$, we define a polyhedron, which is the intersection of all the half-spaces among $H^{\pm}_i$, $i=1,...,s$, that contain $\bx$. Clearly, all such polyhedra define a finite polyhedral structure on $\R^{2n}$, and $Y^{2n-1}$ is a union of entire cells.

Denote by $X^{2n-1}$ the union (maybe empty) of open $(2n-1)$-dimensional cells of $Y^{2n-1}$, and then define $X^{2n-2}:=
Y^{2n-1}\setminus X^{2n-1}$, which is a finite polyhedral complex of dimension $\le2n-2$.

\begin{lemma}\label{cl2} Under the hypotheses of Lemma \ref{cl3}, suppose that $X^{2n-1}\ne\emptyset$.
Then, for each $\bx\in X^{2n-1}$, the preimage $(\Ev^e)^{-1}(\bx)$ consists of regular classes,
or classes $[(\Gamma,h,\bp)]$ such that
\begin{enumerate}\item[(a)] $(\Gamma,h,\bp)$ is of V-type $(\om',\on')$, where $m'_i=m_i$ for all $i\ge1$ except for $m'_{i_1}=m_{i_1}-1$, and
$n'_i=n_i$ for all $i\ge 0$ except for $n'_{i_2}=n_{i_2}-1$ and $n'_{i_1+i_2}=n_{i_1+i_2}+1$;
furthermore, exactly one
connected component of $\Gamma\setminus\bp$ is not regular;
\item[(b)] $\Gamma$ is of V-type $(\om',\on)$, where $m'_i=m_i$ for all $i\ge1$ except either for $m'_{i_1}=m_{i_1}-2$, $m'_{2i_1-2}=m_{2i_1-2}+1$ with some $i_1\ge1$, or for $m'_{i_1}=m_{i_1}-1$, $m'_{i_2}=m_{i_2}-1$, $m'_{i_1+i_2-2}=m_{i_1+i_2-2}+1$ with some $i_2>i_1\ge1$;
furthermore, exactly one connected component of $\Gamma\setminus\bp$ is not regular.
\end{enumerate}
\end{lemma}

{\bf Proof.} By construction, a non-regular element $[(\Gamma,h,\bp)]\in(\Ev^e)^{-1}(\bx)$ is a limit
of regular classes $[(\Gamma_t,h_t,\bp_t)]$, $0<t<\eps$, and is obtained by vanishing of exactly of the parameters in the corresponding cell of
${\mathcal M}^{\;e}_{0,\on,\om}(\R^2,\Delta)$. If the vanishing parameter is the length of a segment joining a marked point $p_{k,t}$ and a vertex $V_t\in\Gamma_t^0
\setminus\bp_t$, then we
get to the case (a). The only non-regular component of $\Gamma\setminus\bp$ is the component which contains the limit of the edge of $\Gamma_t\setminus\bp_t$, which
is not incident to $p_{k,t}$ and is regularly oriented towards $V_t$. If the vanishing parameter is the length of the edge joining two vertices of $\Gamma_t^0\setminus\bp_t$ and
not containing points of $\bp_t$, then we get to the case (b). The only non-regular component of $\Gamma\setminus\bp$ is that with the vertex appeared in the collision of two vertices of $\Gamma_t^0\setminus\bp_t$: the regularity fails, since the new vertex is incident to three incoming edges. Note that no two of these three edges have collinear
directing vectors, since otherwise the dimension of the corresponding cell of $Y^{2n-1}$ would not exceed $2n-2$.
\proofend

By $\widehat{\mathcal M}^{\;e,\lab}_{0,\on,\om}(\R^2,\Del^\lab)$ we denote
the moduli space of labeled $n$-marked plane rational
tropical curves that project to $\widehat{\mathcal M}^{\;e}_{0,\on,\om}(\R^2,\Delta)$.

\subsection{Refined multiplicity of a regular plane rational marked tropical curve}\label{csec2}
Throughout this section, we fix a standard basis in $\R^2$, and
for any $\ba=(a_1,a_2)$, $\bb=(b_1,b_2)\in\R^2$,
set $\ba\wedge\bb=\det\left(\begin{matrix}a_1&a_2\\ b_1&b_2\end{matrix}
\right)$.
For any $\alpha\in\R$ and a formal parameter $y$, define
\begin{equation}[\alpha]_y^-=\frac{y^{\alpha/2}-y^{-\alpha/2}}{y^{1/2}-y^{-1/2}},\quad
[\alpha]_y^+=\frac{y^{\alpha/2}+y^{-\alpha/2}}{y^{1/2}+y^{-1/2}}\ .
\label{ce7}\end{equation}

Now we introduce an additional restriction:
\begin{equation}\sum_{i\ge2}m_i\le1\ .\label{ce20}\end{equation}
It means that all but at most one unmarked vertices are trivalent.

\begin{remark}\label{cr3}
The refined multiplicity of plane tropical curves which we give below naturally extends to arbitrary $\om$ and $\on$ satisfying (\ref{ce2}).
However the invariance statement holds only under restriction (\ref{ce20}). We do not know how to correct the refined multiplicity in order to
obtain an invariant in the general case.
\end{remark}

 Let $[(\Gamma,h,\bp)]\in
{\mathcal M}^{\;e}_{0,\on,\om}(\R^2,\Delta)$, and let
$(\Gamma^{\lab},h,\bp)$ be one of the labelings of
$(\Gamma,h,\bp)$.
We start with defining a refined cuspidal multiplicity $RCM_y(\Gamma,h,\bp,V)$ (depending on
a formal parameter $y$) for each vertex $V\in\Gamma^0$.

\smallskip

\noindent {\it (1) Refined cuspidal multiplicity of a trivalent unmarked vertex.}
Suppose that $V\in\Gamma^0$ is trivalent and the regularly oriented edges $E_1,E_2\in\Gamma^1$ incident
to $V$ are incoming. Define the {\it Mikhalkin's multiplicity} of the vertex $V$ by
(cf. \cite[Definition 2.16]{Mi})
$$\mu(\Gamma,h,\bp,V)=|\ba_1\wedge\ba_2|,\quad\text{where}\quad\ba_i=D(h\big|_{E_i})(\ba_V(E_i)),\ i=1,2\ .$$
Following \cite{BG}, we put
\begin{equation}
RCM_y(\Gamma,h,\bp,V)=
[\mu(\Gamma,h,\bp,V)]_y^-\ .
\label{epsi4}\end{equation}

\smallskip

{\it (2) The function $\mu^+_y(A)$.} We recall here the definition of the function
$\mu^+_y(A)$ for any balanced sequence $A=(\ba_i)_{i=1,...,r}$, $r\ge2$, $\ba_i\in\R^2$, $i=1,...,r$, as given in
\cite[Section 2.1, item (2)]{BS}.
If $r=2$, we set $\mu^+_y(A)=1$. If $r=3$, we set
$\mu^+_y(A)=[|\ba_1\wedge\ba_2|]^+_y$. Note that, due to the balancing condition, this definition does not depend
on the choice of the order in the sequence $A$.
If $r\ge4$, then, for each pair $1\le i<j\le m$, we form the two balanced sequences
\begin{itemize}\item $A'_{ij}$ consisting of the vectors $\ba_k$, $1\le k\le r$, $k\ne i,j$, and one more vector
$\ba_{ij}:=\ba_i+\ba_j$,
\item $A''_{ij}=(\ba_i,\ba_i,-\ba_{ij})$.
\end{itemize} Then we set
\begin{equation}\mu^+_y(A)=\sum_{1\le i<j\le m}\mu^+_y(A'_{ij})\cdot\mu^+_y(A''_{ij})\ .\label{ce8}\end{equation}
It is easy to see that $\mu^+_y(A)$ does not depend on the choice of the order in $A$.

\smallskip

{\it (3) The refined cuspidal multiplicity of a marked vertex.} Given a marked vertex $V\in\Gamma^0\cap\bp$
and the directing vectors
$\ba_1,...,\ba_r$ of all the edges incident to, we set (cf. \cite[Formula (9)]{BS})
\begin{equation}RCM_y(\Gamma,h,\bp,V)=\mu^+_y(A_V),\quad A_v=(\ba_1,...,\ba_r)\ .\label{ce21}\end{equation}

\smallskip
{\it(4) The refined cuspidal multiplicity of an unmarked vertex of valency $\ge4$.} Given an unmarked vertex $V\in\Gamma^0\setminus\bp$ of
valency $r\ge4$ and the directing vectors $\ba_1,...,\ba_r$ of all its incident edges of $\Gamma$, ordered so that $\ba_1,\ba_2$ direct the edges
regularly oriented towards $V$, we set
\begin{equation}RCM_y(\Gamma,h,\bp,V)=[|\ba_1\wedge\ba_2|]^-_y\cdot\mu^+_y(A'_V),\quad A'_V=(\ba_1+\ba_2,\ba_3,...,\ba_r)\ .\label{ce9}\end{equation}

\smallskip
{\it (5) The refined cuspidal multiplicity of a regular plane rational marked tropical curve.}
Given $[(\Gamma,h,\bp)]\in{\mathcal M}^e_{0,\on,\om}(\R^2,\Delta)$, define
\begin{equation}
RCM_y(\Gamma^{\lab},h,\bp)=
\prod_{V\in\Gamma^0}RCM_y(\Gamma,h,\bp,V),
\quad RCM_y(\Gamma,h,\bp)=\frac{RCM_y(\Gamma^{\lab},h,\bp)}{|\Aut(\Gamma,h,\bp)|}\ .\label{ce10}\end{equation}

\subsection{The invariance statement}\label{csec4}

\begin{theorem}\label{ct1}
Let $\Delta\subset\Z^2\setminus\{0\}$ be a balanced, nondegenerate multiset, $\om\in\Z^{\infty,*}_+$, $\on\in\Z^\infty_+$, and let (\ref{ce3}),
(\ref{ce2}) and the restriction (R) hold. Then the expression
\begin{equation}RC_y(\Delta,\om,\bx):=\sum_{[(\Gamma,h,\bp)]\in(\Ev^e)^{-1}(\bx)}RCM_y(\Gamma,h,\bp)
\label{ce13}\end{equation}
does not depend on the choice of $\bx\in\R^{2n}\setminus Y^{2n-1}$ with $Y^{2n-1}$ defined by (\ref{ce11}).
\end{theorem}

\begin{remark}\label{cr1}
In case of $\Delta$ primitive, $n_i=0$ for all $i\ge1$, and $m_i=0$ for all $i\ge2$ (i.e., $(\Gamma,h,\bp)$ trivalent without marked vertices), the cuspidal invariant $RC_y(\Delta,\on,\om)$ coincides with the Block-G\"ottsche refined invariant
$N^{\Delta,\delta}_{trop}(y)$ for $\delta$ chosen so that the counted tropical curves are rational \cite{BG}.

In case of $\Delta$ primitive, $n_i=m_i=0$ for all $i\ge2$ (i.e., $(\Gamma,h,\bp)$ trivalent but with some vertices marked), the
invariant $RC_y(\Delta,\on,\om)$ coincides with refined broccoli invariant as defined by
G\"ottache and Schroeter \cite{GS}.

At last, in case of $m_i=0$ for all $i\ge2$ (i.e., all unmarked vertices trivalent), the invariant $RC_y(\Delta,\on,\om)$
coincides with the refined descendant invariant defined in \cite{BS}. The only novelty of the present note is that we allow one unmarked vertex
of arbitrary valency.
\end{remark}

Similarly to \cite[Proposition 2.4]{BS}, our invariant $RC_y(\Delta,\om)$ is often a rational function of $y$:

\begin{proposition}\label{cp1}
If under hypotheses of Theorem \ref{ct1}, in addition, $\Delta\subset\Z^2\setminus2\Z^2$ (i.e., does not contain even vectors),
then we have
\begin{equation}RC_y(\Delta,\om)=\frac{F(y+y^{-1})}{(y+2+y^{-1})^k}\ ,\label{epsi20}\end{equation}
where $k\ge0$ and $F$ is a nonzero polynomial of degree
$$\deg F=|\Int P(\Delta)\cap\Z^2|+\frac{|\partial P(\Delta)\cap\Z^2|-|\Delta|}{2}+k\ ,$$ where $P(\Delta)$ is the Newton polygon constructed in the proof of
Lemma \ref{cl3}. Furthermore,
\begin{equation}k\le\sum_{i\ge2}i(n_{2i}+n_{2i+1})+\frac{1}{2}\sum_{j\ge4}(j-3)m_j\ .
\label{ce12}\end{equation}
\end{proposition}

{\bf Proof.} The argument used in the proof of \cite[Proposition 2.4]{BS} word-for-word applies in the considered situation. We only make
a couple of comments.
The computation of $\deg F$ uses the construction of a regular tropical curve in
the proof of Lemma \ref{cl3} and also Remark \ref{cr2}. The last summand in the right-hand side of
(\ref{ce12}) (as compared with \cite[Inequality (14)]{BS}) comes from the fact that an unmarked vertex of
valency $j>3$ contributes to the denominator at most $j-3$ factors $y^{1/2}+y^{-1/2}$. \proofend

In general, the denominator in formula (\ref{epsi20}) is unavoidable as noticed in \cite[Corollary 3.3]{BS}.

\subsection{Proof of the invariance}\label{csec5}
It will be convenient to consider labeled tropical curves.
In view of formulas (\ref{ce5}) and (\ref{ce10}), the invariance of $RC_y(\Delta,\on,\om,\bx)$ is equivalent to the
invariance of $RC^{\lab}_y(\Delta,\on,\om,\bx)$. \footnote{The latter expression is defined by formula (\ref{ce13}), where we sum up over all labeled corves.}

So, we choose two generic configurations $\bx(0),\bx(1)\in\R^{2n}\setminus Y^{2n-1}$. There exists a continuous path $\bx(t)\in\R^{2n}$, $0\le t\le 1$, connecting the chosen configurations, that avoids
$X^{2n-2}$, but may finitely many times hit cells of
$X^{2n-1}$, which may cause changes in the structure of $(\Ev^e)^{-1}(\bx(t))$. We shall consider all possible wall-crossing phenomena
and verify the constancy of
$RC^{\lab}_y(\Delta,\om,\om,\bx(t))$ (as a function of $t$) in these events.

To relax notations we simply denote labeled tropical curves by $(\Gamma,h,\bp)$
and write $RC^{\lab}_y(t)$ for $RC^{\lab}_y(\Delta,\om,\bx(t))$.

Let $\bx(t^*)$ be generic in an $(2n-1)$-dimensional cell of $X^{2n-1}$. Denote by $H_0$
the germ of this cell at $\bx(t^*)$ and by $H_+,H_-\subset\R^{2n}$ the germs of the halfspaces
with common boundary $H_0$.
Let $T^*=(\Gamma,h,\bp)\in(\Ev^e)^{-1}(\bx(t^*))$
be as described in Lemma \ref{cl2}(a,b), and let $F_0\subset\widehat{\mathcal M}^{\;e,\lab}_{0,\om}(\R^2,\Delta)$
be the germ at $T^*$ of the $(2n-1)$-cell projected by $\Ev^e$
onto $H_0$. We shall analyze the $2n$-cells of $\widehat{\mathcal M}^{\;e,\lab}_{0,\om}(\R^2,\Delta)$
attached to $F_0$, their projections onto $H_+,H_-$, and prove the
constancy of $RC^{\lab}_y(t)$, $t\in(t^*-\eta,t^*+\eta)$, $0<\eta\ll1$.

\smallskip

{\bf (1)} Suppose that $T^*$ is as in Lemma \ref{cl2}(a), i.e., it has a marked point $p_1$ at a vertex $V\in\Gamma^0$ of valency $i_1+i_2+2$, with
incident edges $E_0,...,E_{i_1+i_2+1}$ directed by the vectors $\ba_j:=\ba_V(E_j)$, $0\le j\le i_1+i_2+1$, and we assume that the limit of the regular orientation is such that $E_0$ is incoming and all other edges incident to $V$ are outgoing. Without loss of generality we can suppose that the path $\bx(t)$, in a neighborhood of $t^*$, is as follows: $x_1=h(p_1)\in\R^2$ moves along a smooth germ transversal to the fixed line $L$ through the segment $h(E_0)$, while $\bx\setminus\{x_1\}$ remains fixed.

Assume that $i_2=0$. Then, in the deformation, the marked point $p_1$ moves from $V$ to one of the edges $E_1,...,E_{i_1+1}$. Note that the sign of $\ba_0\wedge\ba_j$ determines whether the tropical curve
with a marked point on $E_j$, $1\le j\le i_1+1$, is mapped to $H_+$ or $H_-$. Hence, in view of the former formula in
(\ref{ce10}) and formula (\ref{ce9}), the constancy of $RC^{\lab}_y(t)$ is equivalent to the relation
\begin{equation}\sum_{j=1}^{i_1+1}[\ba_0\wedge\ba_j]^-_y\cdot\mu_y^+(A_j)=0,\quad \text{where}\ A_j=(\ba_0+\ba_j,(\ba_k)_{k\ne0,j})\ .
\label{ce14}\end{equation} If $i_1=1$, the balancing condition, which reads $\ba_0+\ba_1+\ba_2=0$, and the definition $\mu^+_y(A_1)=\mu^+_y(A_2)=1$ imply
(\ref{ce14}). If $i_1\ge2$, then (\ref{ce14}) is equivalent to \cite[Formula (18)]{BS}.

\begin{figure}
\setlength{\unitlength}{1.0mm}
\begin{picture}(110,60)(-15,0)
\thicklines

\put(10,45){\line(1,1){5}}\put(10,50){\line(1,0){5}}
\put(10,55){\line(1,-1){5}}\put(15,50){\vector(1,0){10}}
\put(25,50){\line(2,1){10}}\put(35,55){\line(1,1){5}}
\put(35,55){\line(1,0){5}}\put(35,55){\line(1,-1){5}}
\put(25,50){\vector(0,-1){5}}

\put(70,45){\line(1,1){5}}\put(70,50){\line(1,0){5}}
\put(70,55){\line(1,-1){5}}\put(75,50){\line(1,0){10}}
\put(85,50){\line(1,1){5}}\put(85,50){\line(1,-1){5}}
\put(85,50){\line(2,-1){10}}
\put(95,45){\vector(0,-1){5}}\put(95,45){\vector(1,0){5}}

\put(15,15){\vector(0,-1){5}}\put(15,15){\vector(-2,1){5}}
\put(15,15){\line(1,0){10}}\put(25,15){\line(2,1){10}}
\put(25,15){\line(1,-1){5}}\put(25,15){\line(1,1){5}}
\put(25,15){\line(-1,-1){5}}\put(25,15){\line(-1,1){5}}
\put(35,20){\vector(1,1){5}}\put(35,20){\vector(1,-1){5}}

\put(75,15){\line(-2,1){5}}\put(75,15){\line(-1,1){5}}
\put(75,15){\line(-1,-1){5}}\put(75,15){\line(1,1){5}}
\put(75,15){\line(1,-1){5}}\put(75,15){\line(1,0){10}}
\put(85,15){\vector(2,1){10}}\put(95,20){\vector(1,1){5}}
\put(95,20){\vector(1,-1){5}}\put(85,15){\vector(0,-1){5}}

\put(22,0){\rm (c)}\put(80,0){\rm (d)}\put(22,35){\rm (a)}
\put(80,35){\rm (b)}
\put(13.5,48.8){$\bullet$}\put(73.5,48.8){$\bullet$}
\put(3,49){$J\ \Bigg\{$}\put(42,54){$\Bigg\}\ I$}\put(17,52){$\ba'_0$}
\put(17,46){$E'_0$}\put(26,44){$\ba_0$}\put(63,49){$J\ \Bigg\{$}
\put(96,39){$\ba_0$}\put(101,44){$\ba_k$}
\put(11,8){$\ba_j$}\put(5,18){$\ba_k$}
\put(41,13){$\ba_{j_1}$}\put(41,24){$\ba_{j_2}$}\put(101,13){$\ba_{j_1}$}\put(101,24){$\ba_{j_2}$}
\put(86,8){$\ba_j$}

\end{picture}
\caption{Geometric illustration to the invariance statement}\label{cfig1}
\end{figure}
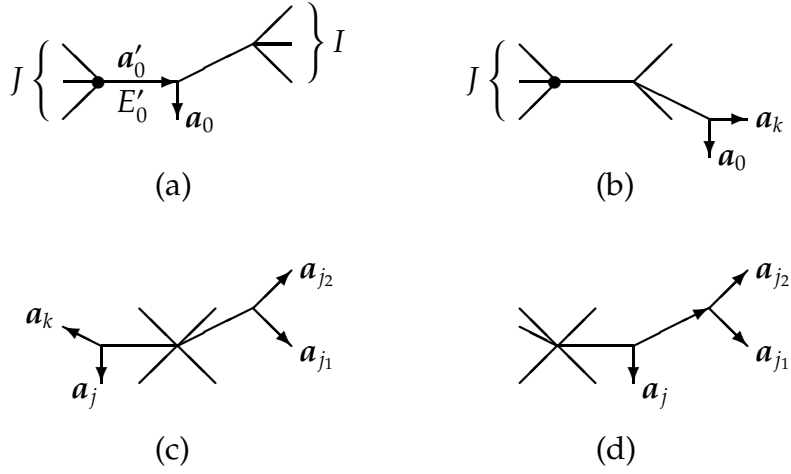

So, assume that $i_2\ge1$. Then, in the deformation, $V$ splits into a marked $(i_2+2)$-valent vertex $p_1$ and an unmarked
$(i_1+2)$-valent vertex $V'$ mapped to the line $L$. Denote by $E'_0$ the edge connecting the vertices $V,V'$ of the deformed curve. The sign of
$\ba_0\wedge\ba'_0$, where $\ba'_0=\ba_V(E'_0)$, determines whether the $\Ev^e$-image of the deformed curve belongs to $H_+$ or to $H_-$. Then the
sought constancy will follow from the relation (see Figure \ref{cfig1}(a))
\begin{equation}\sum_{\renewcommand{\arraystretch}{0.6}
\begin{array}{c}
\scriptstyle{I\cup J=\{1,...,i_1+i_2+1\}}\\
\scriptstyle{|I|=i_1,|J|=i_2+1}\end{array}}[\ba_0\wedge\ba'_0]^-_y\cdot\mu^+_y(A_I)\cdot\mu^+_y(B_J)=0\ ,\label{ce22}\end{equation}
where
$$A_I=(\ba_0,-\ba'_0,(\ba_s\ :\ s\in I)),\quad B_J=(\ba'_0,(\ba_s\ :\ s\in J)),\quad\ba'_0=-\sum_{s\in J}\ba_s\ .$$
Using \cite[Formula (18)]{BS}, we rewrite (\ref{ce22}) in the form (see Figure \ref{cfig1}(b))
$$\sum_{\renewcommand{\arraystretch}{0.6}
\begin{array}{c}
\scriptstyle{I\cup J=\{1,...,i_1+i_2+1\}}\\
\scriptstyle{|I|=i_1,|J|=i_2+1}\end{array}}\sum_{k\in I}[\ba_0\wedge\ba_s]^-_y\cdot\mu^+(A_{I,k})\cdot\mu^+_y(B_J)=0,\quad
A_{I,k}=(\ba_0+\ba_k,-\oa'_0,(\ba_s\ :\ s\in I\setminus\{k\}))\ ,$$ or, equivalently, as
\begin{equation}\sum_{k=1}^{i_1+i_2+1}\left([\ba_0\wedge\ba_k]^-_y\cdot
\sum_{\renewcommand{\arraystretch}{0.6}
\begin{array}{c}
\scriptstyle{I\cup J=\{1,...,i_1+i_2+1\}}\\
\scriptstyle{k\in I,|I|=i_1,|J|=i_2+1}\end{array}}\mu^+_y(A_{I,k})\mu^+_y(B_J)\right)=0\ .\label{ce23}\end{equation}
For a given $k$, the term $\sum_{I,J}\mu^+_y(A_{I,k})\mu^+_y(B_J)$ in the left-hand side of (\ref{ce23}) can be written
(cf. \cite[Section 2.5, proof of Lemma 2.5]{BS}) as the sum of the expressions $\mu^+_{y,\alpha}(C_k,\ba_0+\ba_k,E)$,
where $C_k=(\ba_0+\ba_k,(\ba_s\ :\ 1\le s\le i_1+i_2+1,\ s\ne k))$, $\alpha$ runs over all combinatorial types of
trivalent trees having $i_1+i_2+1$ leaves and containing a point, whose complement consists of two trees with $i_1+1$ and $i_2+2$ leaves, respectively,
which $E$ runs over the leaves of the former subtree. It follows from \cite[Formula (25)]{BS} that
\begin{equation}\sum_{\renewcommand{\arraystretch}{0.6}
\begin{array}{c}
\scriptstyle{I\cup J=\{1,...,i_1+i_2+1\}}\\
\scriptstyle{k\in I,|I|=i_1,|J|=i_2+1}\end{array}}\mu^+_y(A_{I,k})\mu^+_y(B_J)=\Phi_1(z)\sum_{\tau\in S_{i_1+i_2}}z^{\tau\Lambda(C_k)}\ ,
\label{ce24}\end{equation}
where $C_k=\{\ba_s\ :\ 1\le s\le i_1+i_2+1,\ s\ne k\}$, $z^2=y$, $S_k$ is the permutation group of $k$ elements, $C=(\bb_s\}_{1\le s\le|C|}$, and
$$\tau\Lambda(C)=\sum_{1\le s<t\le|C|}\bb_{\tau(s)}\wedge\bb_{\tau(t)}\ .$$
Plugging (\ref{ce24}) to (\ref{ce23}) and using relations
$$[\ba_0\wedge\ba_k]^-_y=\frac{z^{\ba_0\wedge\ba_k}-z^{\ba_k\wedge\ba_0}}{z-z^{-1}},
\quad\ba_0=-\ba_1-...-\ba_{i_1+i_2+1}\ ,$$ we obtain in the left-hand side of (\ref{ce23})
$$\Phi_2(z)\left(\sum_{k=1}^{i_1+i_2+1}\sum_{\tau\in S_{i_1+i_2}}z^{\tau\Lambda(C_k)+\sum_{s\ne0}\ba_k\wedge\ba_s}
-\sum_{k=1}^{i_1+i_2+1}\sum_{\tau\in S_{i_1+i_2}}z^{\tau\Lambda(C_k)+\sum_{s\ne0}\ba_s\wedge\ba_k}\right)$$
$$=\Phi_2(z)\left(\sum_{\sigma\in S_{i_1+i_2+1}}z^{\sigma\Lambda(C)}-\sum_{\sigma\in S_{i_1+i_2+1}}z^{\sigma\Lambda(C)}\right)=0$$
(where $C=\{\ba_1,...,\ba_{i_1+i_2+1}\}$).

\smallskip
{\bf (2)} Suppose that $T^*=(\Gamma,h,\bp)$ is as in Lemma \ref{cl2}(b), i.e., it results from a collision of two unmarked vertices of
valency $3$ and $r\ge3$.
Then $\Gamma$ has an unmarked vertex $V$ of valency $r+1$.
Let $E_j\in\Gamma^1$, $j=1,...,r+1$, be the edges incident to $V$, and the limit of the regular orientation is such that
$E_1,E_2,E_3$ are incoming, while the other edges are outgoing. Denote $\ba_j=\ba_V(E_j)$, $j=1,...,r+1$. We use the same symbols $E_j$
for the corresponding edges of curves $T(t)\in{\mathcal M}^{\;e,\lab}_{0,\on,\om}(\R^2,\Delta)$ obtained in a deformation of $T^*$ along the path
$\bx(t)$, $t\in(t^*-\eta,t^*+\eta)$, no confusion will arise.

The list of possible curves
$T(t)$ is as follows:
\begin{itemize}\item either, for some $1\le j\le3$,
a curve $T(t)$ has a trivalent vertex $V_1$ incident to the edges $E_j$, $E_k$, $k=4,...,r+1$, and the edge $E_0$ that joins $V_1$
with the vertex $V_2$; in turn, $V_2$ is incident to $E_0$, $E_{j_1}$, $E_{j_2}$, where $\{1,2,3\}\setminus\{j\}=\{j_1,j_2\}$, and $E_s$, $s=
4,...,r+1$, $s\ne k$;
\item or, for some $1\le j\le 3$,
a curve $T(t)$ has a trivalent vertex $V_1$ incident to the edges $E_{j_1}$, $E_{j_2}$, and the edge $E_0$ that joins $V_1$
with the vertex $V_2$; in turn, $V_2$ is incident to $E_0$, $E_s$, $s=4,...,r+1$.
\end{itemize}
The regular orientation of $E_0$ is given in the former and in the latter case by the vectors
$$\ba_0=\ba_{V_2}(E_0)=\ba_j+\ba_k\quad\text{and}\quad\ba_0=\ba_{V_1}(E_0)=\ba_{j_1}+\ba_{j_2}\ ,$$
respectively. In both the cases, the sign of $\ba_j\wedge\ba_0$ determines whether $\Ev^e(T(t))$ belongs to $H_+$ or $H_-$.
Introduce $\eps_j=\pm1$, $j=1,2,3$, so that $\eps_j\cdot\sign(\ba_j\wedge\ba_0)=1$ points to $H_+$ for all $j=1,2,3$.
Then the required constancy relation reads
\begin{equation}\sum_{j=1}^3\eps_j\cdot[|\ba_{j_1}\wedge\ba_{j_2}|]^-_y\cdot
\left(\sum_{k=4}^{r+1}[\ba_j\wedge\ba_k]^-_y\cdot\mu^+_y(A_k)+
[\ba_j\wedge(-\ba_{j_1}-\ba_{j_2})]^-_y\cdot\mu^+_y(A)\right)\ ,\label{ce4}\end{equation}
where
$$A_k=(\ba_{j_1}+\ba_{j_2},\ba_j+\ba_k,(\ba_l\ :\ l\in K\setminus\{k\})),\quad
A=(\ba_1+\ba_2+\ba_3,(\ba_l)_{l\in K})$$
(see Figure \ref{cfig1}(c,d)).

If $r=3$, then (\ref{ce4}) turns into
\begin{equation}\sum_{j=1}^3\eps_j\cdot[\oa_j\wedge(-\oa_{j_1}-\oa_{j_2})]^-_y
\cdot[|\oa_{j_1}\wedge\oa_{j_2}|]^-_y=0\ ,
\label{ce6}\end{equation}
which reflects a collision of two trivalent vertices with Block-G\"ottsche refined multiplicities, and in which case (\ref{ce6}) appears to be a particular case of the invariance statement in \cite[Theorem 1]{IM} (see a detailed treatment in \cite[Pages 5313-5316]{IM}).

If $r\ge4$, relation \cite[Formula (18)]{BS} (cf. also (\ref{ce14})) yields that
$$\sum_{k=4}^{r+1}\left([\ba_j\wedge\ba_k]^-_y\cdot\mu^+_y(A_k)\right)=-[\ba_j
\wedge(\ba_{j_1}+\ba_{j_2})]^-_y\cdot\mu^+_y(A)\ ,$$
we obtain in the left-hand side of (\ref{ce4})
$$2\sum_{j=1}^3\left(\eps_j\cdot[|\ba_{j_1}\wedge\ba_{j_2}|]^-_y\cdot[\ba_j
\wedge(-\ba_{j_1}-\ba_{j_2})]^-_y\right)\ ,$$ which vanishes in view of (\ref{ce6}).

\section{On refined count of rational tropical curves in $\R^m$}

From now on we suppose that $m\ge3$, that $\Del\subset\Z^m\setminus\{0\}$ is a balanced, non-degenerate multiset, and that
\begin{equation}n=|\Del|-1\ .\label{se1}\end{equation}

{\bf(1)} Introduce the evaluation maps
$$\Ev:{\mathcal M}_{0,n}(\R^m,\Del)\to\R^{mn},\quad \Ev[(\Gamma,h,\bp)]=h(\bp)\ ,$$
$$\Ev:{\mathcal M}^\lab_{0,n}(\R^m,\Del^\lab)\to\R^{mn},\quad \Ev[(\Gamma^\lab,h,\bp)]=h(\bp)\ .$$
By Lemma \ref{sl1}, $\dim{\mathcal M}_{0,n}(\R^m,\Del)=|\Del|+m+n-3=2n+m-2$. We call an open $(2n+m-2)$-dimensional cell of
${\mathcal M}_{0,n}(\R^m,\Del)$, or ${\mathcal M}^\lab_{0,n}(\R^m,\Del^\lab)$ {\it enumeratively essential}, if $\Ev$ injectively takes it to $\R^{mn}$. Denote by
${\mathcal M}^{\;e}_{0,n}(\R^m,\Del)$, resp. ${\mathcal M}^{\;e,\lab}_{0,n}(\R^m,\Del^\lab)$ the union of all enumeratively essential cells of ${\mathcal M}_{0,n}(\R^m,\Del)$, resp. ${\mathcal M}^\lab_{0,n}(\R^m,\Del^\lab)$, and by $\widehat{\mathcal M}^{\;e}_{0,n}(\R^m,\Del)$, resp.
$\widehat{\mathcal M}^{\;e,\lab}_{0,n}(\R^m,\Del^\lab)$ its closure in ${\mathcal M}_{0,n}(\R^m,\Del)$, resp. ${\mathcal M}^\lab_{0,n}(\R^m,\Del^\lab)$.

By \cite[Construction 2.24 and Proposition 7.4]{GKM} (see also \cite[Construction 4.2]{AR}), the push-forward $\Ev_*{\mathcal M}^\lab_{0,n}(\R^m,\Del^\lab)\subset\R^{mn}$ is a tropical variety (cycle) of dimension $(2n+m-2)$. Moreover, since it is of pure dimension, the support $|\Ev_*{\mathcal M}^\lab_{0,n}(\R^m,\Del^\lab)|\subset\R^{mn}$ coincides with
$\Ev(\widehat{\mathcal M}^{\;e,\lab}_{0,n}(\R^m,\Del^\lab))=\Ev(\widehat{\mathcal M}^{\;e}_{0,n}(\R^m,\Del))$, and hence
we can speak of the tropical variety $\Ev_*\widehat{\mathcal M}^{\;e}_{0,n}(\R^m,\Del)$.

We call a class $[(\Gamma,h,\bp)]\in{\mathcal M}^{\;e}_{0,n}(\R^m,\Del)$, where $n$ satisfies (\ref{se1}), {\it regular}, if
$\Gamma$ is trivalent, $\Gamma^0\cap\bp=\emptyset$, at each vertex $V\in\Gamma^0$ the directing vectors of the incident edges
span a two-dimensional subspace of $\R^m$, and
each connected component $K$ of $\Gamma\setminus\bp$ is unbounded. Notice that the regularity condition together with (\ref{se1})
yield that each connected component of $\Gamma\setminus\bp$ contains exactly one unbounded edge (see \cite[Lemma 4.20]{Mi}).

\begin{lemma}\label{sl3}
Every regular class $[(\Gamma,h,\bp)]\in{\mathcal M}_{0,n}(\R^m,\Del)$ belongs to
${\mathcal M}^{\;e}_{0,n}(\R^m,\Del)$.
\end{lemma}

{\bf Proof.}
We can linearly order the components of $\Gamma\setminus\bp$ and reduce the required statement to the following one: Given a regular rational curve $(\Gamma',h',\bp')$ with $s$ ends and $s-1$ marked points on all but one of the ends, and the $h'$-image of one of the marked points is fixed, then the evaluation image
of the remaining $s-2$ marked points fills a $2(s-2)$-dimensional polyhedron as the curve varies in its combinatorial type. This holds, since we have $2(s-2)$ independent parameters: the lengths of $s-3$ finite edges, and
$s-1$ distances from the marked points to their nearest vertices.
\proofend

%By \cite[Proposition 2.1]{NS} there are finitely many types of rational tropical curves of degree $\Delta$ in $\R^m$, and hence the set $A(\Delta)\subset\Z^m
%\setminus\{0\}$ of
%all possible directing vectors of their edges is finite. Pick an $(m-2)$-dimensional subspace $L\subset\R^m$, which is not parallel to any vector in $A(\Delta)$
%and take a connected neighborhood $U(L)\subset Gr(m,m-2)$ such that none of the spaces $L'\in U(L)$ is parallel to a vector in $A(\Delta)$.
Fix a sequence of $(m-2)$-dimensional linear subspaces $L_1,...,L_{n-1}\in U(L)$ defined over $\Q$, and denote by ${\mathcal L}_i$ the family of $(m-2)$-dimensional affine subspaces of $\R^m$ parallel to $L_i$, $i=1,...,n-1$.

For any point $\bx_0\in\R^m$ and a sequence $\bL=(L'_i\in{\mathcal L}_i,\ i=1,...,n-1)$, introduce the set
$${\mathcal M}^{\;e}_{0,n}(\R^m,\Del,\bx_0,\bL)=\left\{[(\Gamma,h,\bp)]\in\widehat{\mathcal M}^{\;e}_{0,n}(\R^m,\Del)\ \Bigg|\
\begin{matrix}
h(p_1)=\bx_0,\\
h(p_{i+1})\in L'_i,\ i=1,...,n-1\end{matrix}\right\}\ .$$

\begin{remark}\label{sr1}
Our constraint is a particular case of an affine constraint considered in \cite[Definition 1.3]{NS}.
%We believe that there exists an invariant which corresponds to a constraint consisting of a point and
%$n-1$ arbitrary tropical varieties (cycles) of codimension $2$ in $\R^m$.
\end{remark}

\begin{lemma}\label{sl2} For any point $\bx_0\in\R^m$ and generically chosen sequence $\bL$ of affine
spaces $L'_i\in{\mathcal L}_i$, $i=1,...,n-1$, the set ${\mathcal M}^{\;e}_{0,n}(\R^m,\Del,\bx_0,\bL)$ is finite and such that
\begin{itemize}\item
each element $[(\Gamma,h,\bp)]\in{\mathcal M}^{\;e}_{0,n}(\R^m,\Del,\bx_0,\bL)$ is regular;
\item
for each point $p_{i+1}\in\bp$, $i=1,...,n-1$, the germ of the line $h(\Gamma,p_{i+1})$ intersects
with $L_i$ only in one point;
\item for an arbitrary small variation $\bL'$ of the configuration $\bx_0,\bL$, the set
${\mathcal M}^{\;e}_{0,n}(\R^m,\Del,x'_0,\bL')$ is in bijection with ${\mathcal M}^{\;e}_{0,n}(\R^m,\Del,\bx_0,\bL)$ so that any element of
the former set is obtained by a small variation of the corresponding element of the latter set within its combinatorial type.
\end{itemize}
\end{lemma}

{\bf Proof.} The closure of $\Ev({\mathcal M}^{\;e}_{0,n}(\R^m,\Del))$ and the affine space $\{\bx_0\}\times\prod_{j=1}^{n-1}L'_j$ are polyhedral complexes of
complementary dimension in $\R^{mn}$. Hence, after a generic shift of $\{\bx_0\}\times\prod_{j=1}^{n-1}L_j$ (which can be realized via shifts of each
factor in $\R^m$) we obtain a finite intersection, whose points are intersections of open cells of top dimension. In particular,
${\mathcal M}^{\;e}_{0,n}(\R^m,\Del,\bx_0,\bL)\subset{\mathcal M}^{\;e}_{0,n}(\R^m,\Del)$, and the two last claims of the lemma follow.

Next, we prove the regularity of the classes
$[(\Gamma,h,\bp)]\in{\mathcal M}^{\;e}_{0,n}(\R^m,\Del,\bx_0,\bL)$. Consider a component $K$ of $\Gamma\setminus\bp$, whose closure $\overline K$ contains $p_1$ (the marked point mapped to $\bx_0$).
It cannot be bounded. Indeed, otherwise, as observed in the preceding paragraph, the restriction $h:K\to\R^3$ is determined up to a finite
choice by the conditions that one univalent vertex of the closure $\overline K$ is mapped to $\bx_0$ and among the
remaining univalent vertices of $\overline K$ all but one are mapped to certain spaces $L_i$ in the sequence $\bL$. However, this means that the aforementioned conditions
a restriction to the position of the space $L_j$ that contains the very last univalent vertex of
$\overline K$ in
contradiction to the pre-assumed general position of $\bL$. The same argument shows that each other component $K'$ of
$\Gamma\setminus\bp$, whose closure shares a point with $\overline K$, cannot be bounded either. Proceeding in this way, we show that all components of $\Gamma\setminus\bp$ are unbounded. Since
$|\pi_0(\Gamma\setminus\bp)|=n+1=|\Del|$, we obtain that each component of $\Gamma\setminus\bp$ contains exactly one unbounded edge.
It remains to show that the images of the edges of $\Gamma$ incident to the same vertex are not collinear. In fact, it is sufficient to
show this for a component $K$ of $\Gamma\setminus\bp$ having a unique (trivalent) vertex. As we noticed above, the $h$-image of one univalent vertex of $\overline K$ is fixed, while the $h$-image of the other one belongs to a tropical line in general position, but the latter requirement failed if
$h(K)$ were inside a straight line.
\proofend

\smallskip

{\bf(2)} By \cite[Proposition 2.1]{NS} there are finitely many types of rational tropical curves of degree $\Delta$ in $\R^m$, and hence the set $A(\Delta)\subset\Z^m
\setminus\{0\}$ of
all possible directing vectors of their edges
is finite. The set of the subspaces $L'\in Gr(m,m-2)$ such that, for any two vectors $\ba,\bb\in A(\Delta)$, we have $\Psi'(\ba)\wedge\Psi'(\bb)\ne0$ as long as
$\ba\wedge\bb\ne0$, where $\Psi':\R^m\to\R^m/L'$ is the projection,
form a complement of a hypersurface in $GR(m,m-2)$. Pick a connected component $U$ of this set and an element $L\in U$.
%Pick an $(m-2)$-dimensional subspace $L\subset\R^m$, which is not parallel to any vector in $A(\Delta)$
%and take a connected neighborhood $U(L)\subset Gr(m,m-2)$ such that none of the spaces $L'\in U(L)$ is parallel to a vector in $A(\Delta)$.

%Let $R\subset\R^m$ be an $m-2$-dimensional linear subspace such that
%$R$ is transversal to any $2$-plane spanned by integral vectors
%(we call such a subspace $\Z$-{\it general}),
Denote by $\Psi:\R^m\to\R^m/L\simeq\R^2$ the projection, and let $\psi\in\Lambda^2\R^2$ be nondegenerate.

Now we define the weight of each class $[(\Gamma,h,\bp)]\in{\mathcal M}^{\;e}_{0,n}(\R^m,\Del,\bx_0,\bL)$ for $\bx_0$ and $\bL$ in general position. By Lemma \ref{sl2}, $[(\Gamma,h,\bp)]$ is regular, and hence the edges of the closure $\overline K$ of any connected component $K$ of $\Gamma\setminus\bp$ can be uniquely (regularly) oriented so that each marked point is a source, the unbounded edge is directed towards infinity, and for each vertex $V\in K\cap\Gamma^0$ exactly two of its incident edges are incoming. For a vertex $V\in\Gamma^0$, let $E_1,E_2$ be the incoming edges of
$\Gamma\setminus\bp$ incident to $V$, ordered so that $\psi(\Psi(\ba_1),\Psi(\ba_2))>0$,
where $\ba_i=\ba_V(E_i)$, $i=1,2$ (the vanishing here is not possible due to the regularity of $[(\Gamma,h,\bp)]$ and the $\Z$-generality of $R$). Set
$$\mu(\Gamma,h,\bp,L;V)=z^{\ba_1\wedge\ba_2}-z^{\ba_2\wedge\ba_1}\ ,$$
$$SI(\Gamma,h,\bp,L)=\prod_{V\in\Gamma^0}\mu(\Gamma,h,\bp,L;V)\ ,$$
and
\begin{equation}SI(\R^m,\Delta,L,\bx_0,\bL)=\sum_{[(\Gamma,h,\bp)]\in{\mathcal M}^{\;e}_{0,n}(\R^m,\Del,\bx_0,\bL)}
SI(\Gamma,h,\bp,L)\ .\label{se2}\end{equation}

\begin{theorem}\label{st1}
%Suppose that $L_1,...,L_{n-1}\in U$. Then
There exists a neighborhood $U(L)$ of $L$ in $U$ such that the expression $SI(\R^m,\Delta,L,\bx_0,\bL)$ does not depend on the choice of
$L_1,...,L_{n-1}\in U(L)$ and of the choice of a generic constraint $\bx_0\in\R^m$,
$\bL\in\prod_{i=1}^{n-1}{\mathcal L}_i$.
\end{theorem}

\begin{remark}\label{sr2}
(1) All the summands in the right-hand side of (\ref{se2})
belong to the ring $\Z[\Lambda^2\Z^m]$. One can obtain more traditional refined invariants via any group homomorphism $\lambda:\Lambda^*\Z^m\to\Z$. that induces a ring homomorphism $\lambda_*:\Z\big[\frac{1}{2}\Lambda^2\Z^m\big]\to\Z[y^{1/2},y^{-1/2}]$.

(2) It would be interesting to understand an enumerative meaning of our invariant.
More precisely, what is the meaning of the limits of $\lambda_*SI(\R^m,\Delta,L)\cdot(z-z^{-1})^{1-n}$ as $z\to1$ or $z\to\sqrt{-1}$ for various homomorphisms $\lambda:\Lambda^2\Z^m\to\Z$ ? Notice that $n-1$ is the number of vertices of each of the curves
$[(\Gamma,h,\bp)]\in{\mathcal M}^{\;e}_{0,n}(\R^m,\Del,\bx_0,\bL)$.
\end{remark}

Note that $SI_y(\R^m,\Delta,L)$ does depend on the choice of $U$, and the invariance may fail for $L$ parallel to a vector in $A(\Delta)$
(we comment on this in the proof of Theorem \ref{st1} below). We suggests also a relaxed refined invariant %which continuously depends on the choice of an $(m-2)$-subspace $L\subset\R^m$ and a nondegenerate form $\psi\in\Lambda^2(\R^m/L)$. Set
by setting
$$\mu^\red_y(\Gamma,h,\bp,L,\psi;V)=[\psi(\Psi(\ba_1),\Psi(\ba_2))]^-_y=
\frac{y^{\psi(\Psi(\ba_1),\Psi(\ba_2))/2}-
y^{\psi(\Psi(\ba_2),\Psi(\ba_1))/2}}
{y^{1/2}-y^{-1/2}}\ ,$$
$$SI^\red_y(\Gamma,h,\bp,L,\psi)=\prod_{V\in\Gamma^0}\mu^\red_y(\Gamma,h,\bp,L,\psi;V)\ ,$$
and
\begin{equation}SI^\red_y(\R^m,\Delta,L,\psi,\bx_0,\bL)=\sum_{[(\Gamma,h,\bp)]\in{\mathcal M}^{\;e}_{0,n}(\R^m,\Del,\bx_0,\bL)}
SI^\red_y(\Gamma,h,\bp,L,\psi)\ .\label{se5}\end{equation}

Thus, we immediately obtain

\begin{corollary}\label{st2}
Under the hypotheses of Theorem \ref{st1}, the expression $SI^\red_y(\R^m,\Delta,L,\psi,\bx_0,\bL)$ does not depend on the choice of a generic constraint $\bx_0\in\R^m$ and
$\bL\in\prod_{i=1}^{n-1}{\mathcal L}_i$.
\end{corollary}

\begin{example}\label{sex1}
Suppose that $L_1=...=L_{n-1}=L$ and that %\{y_1=y_2=0\}$ ($y_1,...,y_m$ being the coordinates in $\R^m$) and that
$\psi(\ba,\bb)=y_1(\ba)y_2(\bb)-y_2(\ba)y_1(\bb)$.
Then $SI^\red_y(\R^m,\Delta,(L_i)_{i=1}^{n-1})$ turns to be the Block-G\"ottsche refined invariant for plane rational curves of degree $\Psi(\Delta)$,
multiplied by $(y^{1/2}-y^{-1/2})^{n-1}$.
\end{example}

\smallskip
{\bf Proof of Theorem \ref{st1}.} We fix the point $\bx_0\in\R^m$, choose two generic configurations $\bL(0),
\bL(1)\in\prod_{i=1}^{n-1}$, join them by a generic path $\{\bL(t)\}_{0\le t\le 1}$, and show that
$SI(\R^m,\Del,L,\bx_0,\bL(t))$ remains constant along this path. More precisely,
for all but finitely many $t\in[0,1]$, the set ${\mathcal M}^{\;e}_{0,n}(\R^m,\Delta,\bx_0,\bL(t))$ satisfies the properties listed in Lemma \ref{sl2}, and hence $SI(\R^m,\Del,L,\bx_0,\bL(t))$ is constant along each connected component of the complement to the above finite exceptional set $F\subset(0,1)$.
%To simplify the analysis of the events corresponding to the points of $F$, we slightly deform the spaces $L_1,...,L_{n-1}$ leaving them rational (i.e., defined over $\Q$), but not containing any of the directing vectors of the edges of the classes $[(\Gamma,h,\bp)]\in\widehat{\mathcal M}^{\;e}_{0,n}(\R^m,\Delta)$.
%This deformation does not affect the claim we want to prove
%(cf. the first paragraph in the proof of Lemma \ref{sl2}) and it ensures that,
Observe that, due to $L_1,...,l_{n-1}\in U$, for $[(\Gamma,h,\bp)]\in{\mathcal M}_{0,n}(\R^m,\Del)$, the intersection of
$h(\Gamma)$ with any of the spaces $L'_i\in{\mathcal L}_i$, $i=1,...,n-1$, is finite.

The points of $F$ correspond to the events when the current constraint $\bL(t)$ meets the $\Ev$-image of a codimension one cell of $\widehat{\mathcal M}^{\;e}_{0,n}(\R^m,\Delta)$ in its generic point. It follows that exactly one of the parameters in the dimension count in the proof of Lemma \ref{sl2} (cf. also \cite[Proof of Theorem 1]{IM}) vanishes. More precisely, for $t^*\in F$, all but one elements of ${\mathcal M}^{\;e}_{0,n}(\R^m,\Del,\bx_0,\bL(t^*))$ are regular, while the remaining element $[(\Gamma,h,\bp)]$ is as follows
\begin{enumerate}\item[(i)] either $\Gamma$ is trivalent, $\bp\cap\Gamma^0=\{p_i\}$ for some $1\le i\le n$, all but one components of $\Gamma\setminus\bp$ are unbounded, while one component is bounded and contains $p_i$ in its boundary;
\item[(ii)]
or $\bp\cap\Gamma=\emptyset$, all but one vertices of $\Gamma$ are trivalent, and one vertex is four-valent and is incident to three incoming edges; all components of $\Gamma\setminus\bp$ are unbounded.
\end{enumerate}
Denote by $\Gamma_V$ the germ of $\Gamma$ at $V$, where $V=p_i$ in case (i), and $V\in\Gamma^0$ is the four-valent vertex in case (ii).

\smallskip

\begin{figure}
\setlength{\unitlength}{1cm}
\begin{picture}(15,3)(0.5,0)
\thinlines

\put(7.2,2){\line(1,1){1}}\put(7.2,2){\line(-1,1){1}}\put(7.2,2){\line(0,-1){1}}
\put(8.8,1.9){$\Longrightarrow$}\put(5,1.9){$\Longleftarrow$}
\put(10.7,1.7){\line(1,1){1}}\put(10.7,1.7){\line(-1,1){0.9}}\put(10.7,1.7){\line(0,-1){0.7}}
\put(3.5,1.7){\line(1,1){0.9}}\put(3.5,1.7){\line(-1,1){1}}\put(3.5,1.7){\line(0,-1){0.7}}
\put(7.1,1.9){$\bullet$}\put(10.9,1.9){$\bullet$}\put(3.1,1.9){$\bullet$}
\put(3.4,0.6){$e_1$}\put(7.1,0.6){$e_1$}\put(10.6,0.6){$e_1$}
\put(2.1,2.8){$e_2$}\put(5.8,2.8){$e_2$}\put(9.3,2.8){$e_2$}
\put(4.1,2.8){$e_3$}\put(7.5,2.8){$e_3$}\put(11.3,2.8){$e_3$}

\end{picture}
\caption{Bifurcation}\label{sf1i}
\end{figure}
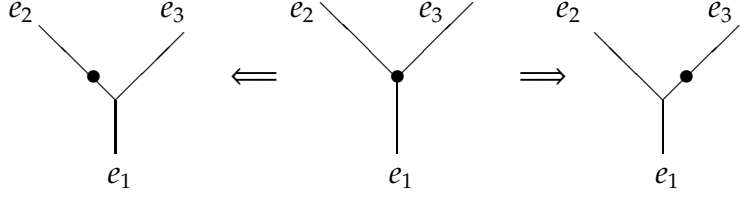

{\it The case (i): Passage through $[(\Gamma,h,\bp)]\in{\mathcal M}^{\;e}_{0,n}(\R^m,\Delta\bx_0,\bL(t^*))$ with $p_i=V\in\Gamma^0$, $1\le i\le n$.} Let $E_1,E_2,E_3\in\Gamma^1$ be incident to $V$. Exactly one of these edges, say, $E_1$, is contained in a bounded component of $\Gamma$. Since $p_i$ cannot move inside $E_1$ by the regularity condition (see Lemma \ref{sl2}), there are two top-dimensional cells of $\widehat{\mathcal M}^{\;e}_{0,n}(\R^m,\Delta)$ attached to the cell containing $[(\Gamma,h,\bp)]$ and corresponding to the moves of $p_i$ into $E_2$, or $E_3$. Since $\Ev_*\widehat{\mathcal M}^{\;e}_{0,n}(\R^m,\Delta)$ is a tropical variety, the balancing condition ensures that the $\Ev$-images of the above top-dimensional cells of $\widehat{\mathcal M}^{\;e}_{0,n}(\R^m,\Delta)$ cover a germ of an $(m+2n-2)$-plane centered at $\Ev[(\Gamma,h,\bp)]$.
Thus, in the plane $\R^2$, which is the target of $\Phi$, the considered passage looks as shown in Figure \ref{sf1i}, where $e_1,e_2,e_3$, and $\widetilde p_i$ denote the $\Psi\circ h$-images of $E_1,E_2,E_3$, and $p_i$, respectively. Assuming that the numbering of $E_1,E_2,E_3$ matches the cyclic order defined by the form $\psi$, we
reduce the constancy of $SI(\R^m,\Del,L,\bx_0,\bL(t))$ for
$t\in(t^*-\eps,t^*+\eps)\setminus\{t^*\}$ to the relation
$$\ba_V(E_1)\wedge\ba_V(E_2)=\ba_V(E_3)\wedge\ba_V(E_1)\ ,$$ which follows from the equality
$\ba_V(E_3)=-\ba_V(E_1)-\ba_V(E_2)$ coming from (\ref{se3}).

\smallskip
{\it The case (ii): Passage through $[(\Gamma,h,\bp)]\in{\mathcal M}^{\;e}_{0,n}(\R^m,\Del,\bx_0,\bL(t^*))$ with
$\Gamma$ having a four-valent vertex.} Let $V\in\Gamma^0$ be a four-valent vertex, $E_1,E_2,E_3,E_4\in\Gamma^1$ the edges incident to $V$, whose $\Psi\circ h$-images $e_1,e_2,e_3,e_4$ (see Figure \ref{sf1}(b,c)) are cyclically ordered so that $\psi(\Psi(\ba_i),\Psi(\ba_{i+1}))>0$, where $\ba_i=\ba_V(E_i)$, $i=1,2,3,4$. We can assume also that, in the limit of the regular orientation, $E_1,E_2,E_3$ are incoming, while $E_4$ is outgoing.

Suppose that no two of
the edges $e_1,e_2,e_3,e_4$ are parallel to each other. In this case, there are three top-dimensional cells of $\widehat{\mathcal M}^{\;e}_{0,n}(\R^m,\Del)$ incident to $[(\Gamma,h,\bp)]$ and corresponding to three ways of splitting of a four-valent vertex into a pair of trivalent vertices (cf. \cite[Page 156]{GM} or \cite[Proof of Theorem 3.7]{GKM}). We can assume that the path $\bL(t)$,
$t\in(t^*-\eps,t^*+\eps)$, is a generic line in $\prod_{i=1}^{n-1}{\mathcal L}_i$. Again, since $\Ev_*\widehat{\mathcal M}^{\;e}_{0,n}(\R^m,\Delta)$ is a tropical variety, the balancing condition ensures that the $\Ev$-images of the above three top-dimensional cells of $\widehat{\mathcal M}^{\;e}_{0,n}(\R^m,\Delta)$ cover the germ of an $(m+2n-2)$-plane $H$ centered at $\Ev[(\Gamma,h,\bp)]$. That is, one of the cells maps to one half-space $H^+$ of $H$, bounded by the considered wall, while the other two cella are mapped to the other half-space $H^-$ of $H$.

Assume for the moment that $L_1=...=L_{n-1}=L$. Under this assumption, all the elements of the constraint are mapped to distinct points in $\R^m/L\simeq\R^2$, and hence we encounter the well-known trifurcation for plane tropical curves matching suitably many point constraints (see, for instance,
\cite[Page 156]{GM} and \cite[Proof of Theorem 3.7]{GKM}). Namely, this trifurcation and the corresponding
subdivisions of the plane quadrangle $Q$, dual to the four-valent vertex, are shown in Figure \ref{sf1}(a). The distribution of the cells of $\widehat{\mathcal M}^{\;e}_{0,n}(\R^m,\Delta)$ mapped to $H^+$ or $H^-$ depends on the position of the position of the unique parallelogram inscribed into $Q$, and it is indicated
in the left part of Figure \ref{sf1}(a). More precisely, this distribution matches the equality
$$\mu_1\mu_2=\mu_3\mu_4+\mu_5\mu_6\ ,$$ where $\mu_1,...,\mu_6$ are twice the Euclidean areas of the triangles as indicated in the right part of Figure \ref{sf1}(a)
(see for details the aforementioned sources). So, depending on the position of the vectors $e_1,...,e_4$ with respect to the inscribed parallelogram, we encounter
the two possible options shown in Figure \ref{sf1}(b,c).
%This yields that, in the deformation induced by this path, the edges $e_1,e_2,e_3$ move in a parallel way with generic constant speeds, while the move of $e_4$ (also in a parallel way) is determined by splitting of $V$ into a pair of trivalent vertices in by the balancing condition for the plane tropical curve
%$(\Gamma,\Psi\circ h)$. A suitable family of coordinate changes in $\R^2$ fixes the lines containing the edges $e_1,e_2$. Then the move of of $e_3$ leads to three types of deformation of the germ of the plane tropical curve $(\Gamma_V,\Psi\circ h)$ as shown in Figure \ref{sf1}(b,c): they are dual to three subdivisions of the Newton quadrangle associated with the four-valent vertex of
%$(\Gamma_V,\Psi\circ h)$, cf. \cite[Figures 6 and 7]{IM} or \cite[Figure 4]{SS} (the difference between Figures \ref{sf1}(b) and (c) is that, in the former case, $e_4$ is dual to a common edge of the quadrangle and the inscribed parallelogram, while, in the latter case, it is not so). The relative position of $e_3$ with respect to $e_1,e_2$ shows that, in both figures, the upper diagram corresponds, say, to $t>t^*$, while the two lower ones to $t<t^*$.
Hence, the constancy of $SI(\R^m,\Del,L,\bx_0,\bL(t))$, $t\in(t^*-\eps,t^*+\eps)$, reduces to the equalities
$$(z^{\ba_1\wedge\ba_2}-z^{\ba_2\wedge\ba_1})(z^{(\ba_1+\ba_2)\wedge\ba_3}-z^{\ba_3\wedge(\ba_1+\ba_2)})$$%\mnote{\ba_2 \to \ba_3}
$$=(z^{\ba_2\wedge\ba_3}-z^{\ba_3\wedge\ba_2})(z^{\ba_1\wedge(\ba_2+\ba_3)}-z^{(\ba_2+\ba_3)\wedge\ba_1})
+(z^{\ba_3\wedge\ba_1}-z^{\ba_1\wedge\ba_3})(z^{\ba_2\wedge(\ba_1+\ba_3)}-z^{(\ba_1+\ba_3)\wedge\ba_2})$$
(Figure \ref{sf1}(b)) and
$$(z^{\ba_2\wedge\ba_3}-z^{\ba_3\wedge\ba_2})(z^{\ba_1\wedge(\ba_2+\ba_3)}-z^{(\ba_2+\ba_3)\wedge\ba_1})$$
$$=(z^{\ba_1\wedge\ba_2}-z^{\ba_2\wedge\ba_1})(z^{(\ba_1+\ba_2)\wedge\ba_3}
-z^{\ba_3\wedge(\ba_1+\ba_2)})
+(z^{\ba_1\wedge\ba_3}-z^{\ba_3\wedge\ba_1})(z^{\ba_2\wedge(\ba_1+\ba_3)}-z^{(\ba_1+\ba_3)\wedge\ba_2})$$
(Figure \ref{sf1}(c)), which are clearly valid. If $L_1,...,L_{n-1}$ vary in a sufficiently small neighborhood $U(L)$ of $L$, then the considered three cells of
$\widehat{\mathcal M}^{\;e}_{0,n}(\R^m,\Delta)$ remain to be mapped to the same half-space $H^+$ or $H^-$ as above, and hence again
the constancy $SI(\R^m,\Del,L,\bx_0,\bL(t))$, $t\in(t^*-\eps,t^*+\eps)$, follows.

Suppose that (at least) two of the edges $e_1,e_2,e_3,e_4$ are parallel, that is, either $e_1\|e_3$, or $e_2\|e_4$. Due to the
choice of $L$, %$\Z$-generality of $R$,
the latter relations yield $\ba_1\|\ba_3$ and $\ba_2\|\ba_4$, respectively. In these cases we encounter only two types of deformations of $[(\Gamma,h,\bp)]$, see Figure \ref{sf1}(d,e). Thus, the constancy of $SI(\R^m,\Del,L,\bx_0,\bL(t))$, $t\in(t^*-\eps,t^*+\eps)$, reduces in both the cases to the equality
$$(z^{\ba_2\wedge\ba_1}-z^{\ba_1\wedge\ba_2})(z^{\ba_3\wedge(\ba_1+\ba_2)}-z^{(\ba_1+\ba_2)\wedge\ba_3})=
(z^{\ba_3\wedge\ba_2}-z^{\ba_2\wedge\ba_3})(z^{(\ba_2+\ba_3)\wedge\ba_1}-z^{\ba_1\wedge(\ba_2+\ba_3)})\ .$$ If $\ba_1\|\ba_3$, then it follows from $\ba_1\wedge\ba_3=0$, and if $\ba_2\|\ba_4$, it follows from $\ba_2\wedge\ba_4=0$ and $\ba_4=-\ba_1-\ba_2-\ba_3$.
\proofend

Notice that if $L\not\in U$, in the last paragraph of the proof one may have $\ba_1\not\parallel\ba_3$, resp. $\ba_2\not\parallel\ba_4$, and the required invariance may fail.

\begin{figure}
\setlength{\unitlength}{1cm}
\begin{picture}(15,19)(0.5,0)
\thinlines

\put(1.5,16.6){\line(1,1){1}}\put(1.5,16.6){\line(1,-1){1}}
\put(1.5,16.6){\line(-1,2){0.5}}\put(1.5,16.6){\line(-1,-2){0.5}}

\put(5.2,18.7){\line(1,-1){0.8}}\put(5.2,18.7){\line(-1,-2){0.4}}
\put(5.2,19.1){\line(1,1){0.8}}\put(5.2,19.1){\line(-1,2){0.4}}
\put(5.2,18.7){\line(0,1){0.4}}

\put(5,16.6){\line(-1,-2){0.4}}\put(5,16.6){\line(-1,2){0.4}}
\put(5.4,16.6){\line(1,1){0.6}}\put(5.4,16.6){\line(1,-1){0.6}}
\put(5,16.6){\line(1,0){0.4}}

\put(5,14.5){\line(-1,-2){0.3}}\put(5,14.9){\line(-1,2){0.3}}
\put(5,14.5){\line(1,1){0.7}}\put(5,14.9){\line(1,-1){0.7}}
\put(5,14.5){\line(0,1){0.4}}

\put(3,16.5){$\Longrightarrow$}
\put(3,17.8){$\Nearrow$}\put(3,15.3){$\Searrow$}

%\put(3.3,0){(a)}\put(10.8,0){(b)}

\put(7,18.7){\line(2,1){2}}\put(7,18.7){\line(2,-1){2}}
\put(9,19.7){\line(1,-1){1}}\put(9,17.7){\line(1,1){1}}
\put(7,18.7){\line(1,0){3}}

\put(7,16.6){\line(2,1){2}}\put(7,16.6){\line(2,-1){2}}
\put(9,17.6){\line(1,-1){1}}\put(9,15.6){\line(1,1){1}}
\put(9,15.6){\line(0,1){2}}

\put(7,14.5){\line(2,1){2}}\put(7,14.5){\line(2,-1){2}}
\put(9,15.5){\line(1,-1){1}}\put(9,13.5){\line(1,1){1}}
\put(9,15.5){\line(-1,-1){1}}\put(9,13.5){\line(-1,1){1}}
\put(7,14.5){\line(1,0){1}}

\put(11.7,16.6){\line(2,1){2}}\put(11.7,16.6){\line(2,-1){2}}
\put(13.7,17.6){\line(1,-1){1}}\put(13.7,15.6){\line(1,1){1}}

\put(10.5,16.5){$\Longleftarrow$}
\put(10.7,17.8){$\Nwarrow$}\put(10.7,15.3){$\Swarrow$}

%\put(0.5,16.5){$F_0$}
\put(4.2,18.6){$H^+$}\put(4.2,16.5){$H^-$}\put(4.2,14.6){$H^-$}

\put(8.6,18.9){$\mu_1$}\put(8.6,18.3){$\mu_2$}
\put(8.1,16.5){$\mu_3$}\put(9.2,16.5){$\mu_4$}
\put(7.5,15.1){$\mu_5$}\put(7.5,13.8){$\mu_6$}

%\put(7.2,15){\line(1,1){1}}\put(7.2,15){\line(-1,1){1}}\put(7.2,15){\line(0,-1){1}}
%\put(8.8,14.9){$\Longrightarrow$}\put(5,14.9){$\Longleftarrow$}
%\put(10.7,14.7){\line(1,1){1}}\put(10.7,14.7){\line(-1,1){0.9}}\put(10.7,14.7){\line(0,-1){0.7}}
%\put(3.5,14.7){\line(1,1){0.9}}\put(3.5,14.7){\line(-1,1){1}}\put(3.5,14.7){\line(0,-1){0.7}}
%\put(7.1,14.9){$\bullet$}\put(10.9,14.9){$\bullet$}\put(3.1,14.9){$\bullet$}
%\put(3.4,13.6){$e_1$}\put(7.1,13.6){$e_1$}\put(10.6,13.6){$e_1$}
%\put(2.1,15.8){$e_2$}\put(5.8,15.8){$e_2$}\put(9.3,15.8){$e_2$}
%\put(4.1,15.8){$e_3$}\put(7.5,15.8){$e_3$}\put(11.3,15.8){$e_3$}

\put(6.1,9.8){$e_4$}\put(6.1,11.9){$e_3$}\put(4.4,9.8){$e_1$}\put(4.4,11.9){$e_2$}
\put(6.1,7.8){$e_4$}\put(6.1,9.3){$e_3$}\put(4.2,7.8){$e_1$}\put(4.2,9.4){$e_2$}
\put(5.8,5.9){$e_4$}\put(5.9,7.2){$e_3$}\put(4.3,5.8){$e_1$}\put(4.3,7.3){$e_2$}
\put(2.2,7.4){$e_4$}\put(2.2,9.7){$e_3$}\put(0.6,7.6){$e_1$}\put(0.6,9.5){$e_2$}

\put(13.1,9.8){$e_3$}\put(13.1,11.9){$e_2$}\put(11.4,9.8){$e_4$}\put(11.4,11.9){$e_1$}
\put(13.1,7.8){$e_3$}\put(13.1,9.3){$e_2$}\put(11.2,7.8){$e_4$}\put(11.2,9.4){$e_1$}
\put(12.8,5.9){$e_3$}\put(12.9,7.2){$e_2$}\put(11.3,5.8){$e_4$}\put(11.3,7.3){$e_1$}
\put(9.2,7.4){$e_3$}\put(9.2,9.7){$e_2$}\put(7.6,7.6){$e_4$}\put(7.6,9.5){$e_1$}

\put(3.3,3.9){$e_3$}\put(4,3.4){$e_4$}\put(2.7,2.9){$e_2$}\put(3.3,1.4){$e_1$}
\put(0.5,4.1){$e_3$}\put(1.5,3.4){$e_4$}\put(0,2.8){$e_2$}\put(0.8,1.4){$e_1$}
\put(6,4.1){$e_3$}\put(6.9,3.6){$e_4$}\put(5.2,2.9){$e_2$}\put(6,1.4){$e_1$}

\put(11.3,3.9){$e_2$}\put(12,3.4){$e_3$}\put(10.7,2.9){$e_1$}\put(11.3,1.4){$e_4$}
\put(8.5,4.1){$e_2$}\put(9.5,3.4){$e_3$}\put(8,2.8){$e_1$}\put(8.8,1.4){$e_4$}
\put(14,4.1){$e_2$}\put(14.9,3.6){$e_3$}\put(13.2,2.9){$e_1$}\put(14,1.4){$e_4$}

\put(1.5,8.6){\line(1,1){1}}\put(1.5,8.6){\line(1,-1){1}}
\put(1.5,8.6){\line(-1,2){0.5}}\put(1.5,8.6){\line(-1,-2){0.5}}

\put(5.2,10.7){\line(1,-1){0.8}}\put(5.2,10.7){\line(-1,-2){0.4}}
\put(5.2,11.1){\line(1,1){0.8}}\put(5.2,11.1){\line(-1,2){0.4}}
\put(5.2,10.7){\line(0,1){0.4}}

\put(5,8.6){\line(-1,-2){0.4}}\put(5,8.6){\line(-1,2){0.4}}
\put(5.4,8.6){\line(1,1){0.6}}\put(5.4,8.6){\line(1,-1){0.6}}
\put(5,8.6){\line(1,0){0.4}}

\put(5,6.5){\line(-1,-2){0.3}}\put(5,6.9){\line(-1,2){0.3}}
\put(5,6.5){\line(1,1){0.7}}\put(5,6.9){\line(1,-1){0.7}}
\put(5,6.5){\line(0,1){0.4}}

\put(3,8.5){$\Longrightarrow$}
\put(3,9.8){$\Nearrow$}\put(3,7.3){$\Searrow$}

\put(7,13){(a)}\put(3.5,5){(b)}\put(10.5,5){(c)}
\put(3.2,0){(d)}\put(11.2,0){(e)}

\put(3.5,2.7){\line(0,1){1}}\put(3.5,2.7){\line(0,-1){1}}
\put(3.5,2.7){\line(1,1){0.6}}\put(3.5,2.7){\line(-1,0){0.7}}
\put(4.4,2.6){$\Longrightarrow$}\put(1.7,2.6){$\Longleftarrow$}
\put(6.2,3.1){\line(0,1){0.8}}\put(6,2.7){\line(0,-1){1}}
\put(6.2,3.1){\line(1,1){0.6}}\put(6,2.7){\line(-1,0){0.7}}
\put(6,2.7){\line(1,2){0.2}}
\put(0.6,3.1){\line(0,1){0.8}}\put(1,2.7){\line(0,-1){1}}
\put(1,2.7){\line(1,1){0.6}}\put(0.6,3.1){\line(-1,0){0.6}}
\put(1,2.7){\line(-1,1){0.4}}

\put(11.5,2.7){\line(0,1){1}}\put(11.5,2.7){\line(0,-1){1}}
\put(11.5,2.7){\line(1,1){0.6}}\put(11.5,2.7){\line(-1,0){0.7}}
\put(12.4,2.6){$\Longrightarrow$}\put(9.7,2.6){$\Longleftarrow$}
\put(14.2,3.1){\line(0,1){0.8}}\put(14,2.7){\line(0,-1){1}}
\put(14.2,3.1){\line(1,1){0.6}}\put(14,2.7){\line(-1,0){0.7}}
\put(14,2.7){\line(1,2){0.2}}
\put(8.6,3.1){\line(0,1){0.8}}\put(9,2.7){\line(0,-1){1}}
\put(9,2.7){\line(1,1){0.6}}\put(8.6,3.1){\line(-1,0){0.6}}
\put(9,2.7){\line(-1,1){0.4}}

\put(8.5,8.6){\line(1,1){1}}\put(8.5,8.6){\line(1,-1){1}}
\put(8.5,8.6){\line(-1,2){0.5}}\put(8.5,8.6){\line(-1,-2){0.5}}

\put(12.2,10.7){\line(1,-1){0.8}}\put(12.2,10.7){\line(-1,-2){0.4}}
\put(12.2,11.1){\line(1,1){0.8}}\put(12.2,11.1){\line(-1,2){0.4}}
\put(12.2,10.7){\line(0,1){0.4}}

\put(12,8.6){\line(-1,-2){0.4}}\put(12,8.6){\line(-1,2){0.4}}
\put(12.4,8.6){\line(1,1){0.6}}\put(12.4,8.6){\line(1,-1){0.6}}
\put(12,8.6){\line(1,0){0.4}}

\put(12,6.5){\line(-1,-2){0.3}}\put(12,6.9){\line(-1,2){0.3}}
\put(12,6.5){\line(1,1){0.7}}\put(12,6.9){\line(1,-1){0.7}}
\put(12,6.5){\line(0,1){0.4}}

\put(10,8.5){$\Longrightarrow$}
\put(10,9.8){$\Nearrow$}\put(10,7.3){$\Searrow$}

\end{picture}
\caption{Trifurcation}\label{sf1}
\end{figure}

%{\bf Proof of Theorem \ref{st2}.}
%We literally follow the proof of Theorem \ref{st1}, commenting only on the last considered case (Figures \ref{sf1}(d,e)): for any $L$, we have either $\Psi(\ba_1)\|\Psi(\ba_3)$, or $\Psi(\ba_2)\|\Psi(\ba_4)$, and hence the desired relation
%$$(y^{\psi(\Psi(\ba_2),\Psi(\ba_1))/2}-y^{\psi(\Psi(\ba_1),\Psi(\ba_2))/2})
%(y^{\psi(\Psi(\ba_3),\Psi(\ba_1+\ba_2))/2}-y^{\psi(\Psi(\ba_1+\ba_2),\Psi(\ba_3))/2})$$
%$$=(y^{\psi(\Psi(\ba_3),\Psi(\ba_2))/2}-y^{\psi(\Psi(\ba_2),\Psi(\ba_3))/2})
%(y^{\psi(\Psi(\ba_2+\ba_3),\Psi(\ba_1))/2}-y^{\psi(\Psi(\ba_1),\Psi(\ba_2+\ba_3))/2})\qquad \text{\proofend}$$

%{\ncsc School of Mathematical Sciences \\[-21pt]
%
%Raymond and Beverly Sackler Faculty of Exact Sciences\\[-21pt]
%
%Tel Aviv University \\[-21pt]
%
%Ramat Aviv, 69978 Tel Aviv, Israel} \\[-21pt]
%
%{\it E-mail address}: {\ntt shustin@post.tau.ac.il}


\begin{thebibliography}{99}
\bibitem{AR} L. Allermann and J. Rau.
First Steps in Tropical Intersection Theory. {\it Math. Z.} {\bf 264} (2010), no. 3, 633--670.
\bibitem{BS} L. Blechman and E. Shustin.
Refined descendant invariants of toric surfaces.
{\it Discr. Comput. Geom.} {\bf 62} (2019), no. 1, 180--208.
\bibitem{BG} F. Block and L. G\"ottsche.
Refined curve counting with tropical geometry.
{\it Compos. Math.} {\bf 152} (2016), no. 1, 115--151.
\bibitem{BG1} F. Block and L. G\"ottsche.
Fock spaces and refined Severi degrees. {\it Int. Math. Res. Not.} <b>2016</b>,  no. 21, 6553--6580.
\bibitem{GSh} Y. Ganor and E. Shustin.
{\it Enumeration of plane unicuspidal curves of any genus via tropical geometry}.
Preprint at arXiv:1807.11443.
\bibitem{GKM} A. Gathmann, M. Kerber, and H. Markwig. Tropical fans and
the moduli spaces of tropical curves. {\it Compos. Math.} {\bf 145} (2009), no.
1, 173--195.
\bibitem{GM} A. Gathmann and H. Markwig. The numbers of
tropical plane curves through points in general position.
{\it J. reine angew. Math.} {\bf 602} (2007), 155--177.
\bibitem{GS} L. G\"ottsche and F. Schroether. Refined broccoli invariants.
{\it J. Alg. Geom.} {\bf 28} (2019), 1--41.
\bibitem{GoS} L. G\"ottsche and V. Shende.
Refined curve counting on
complex surfaces. {\it Geom. Topol.} {\bf 18} (2014) 2245--2307.
\bibitem{IM} I. Itenberg and G. Mikhalkin. On Block-G\"ottsche
multiplicities for planar tropical curves. {\it IMRN} {\bf 23} (2013), 5289--5320.
\bibitem{Ma} T. Mandel. {\it Scattering diagrams, theta functions,
and refined tropical curve counts}.
Preprint at arXiv:1503.06183.
\bibitem{MaRu} T. Mandel and H. Ruddat.
{\it Descendant log Gromov-Witten invariants for toric
varieties and tropical  curves}. Preprint at arXiv:1612.02402.
\bibitem{MaRu1} T. Mandel and H. Ruddat.
{\it Tropical quantum field theory, mirror polyvector fields,
and multiplicities of tropical curves}. Preprint at arXiv:1902.07183.
\bibitem{Mi} G. Mikhalkin. Enumerative tropical algebraic geometry in $\R^2$. {\it
J. Amer. Math. Soc.} {\bf 18} (2005), 313--377.
\bibitem{Mi3} G. Mikhalkin. Tropical Geometry and its applications.
Sanz-Solé, Marta (ed.) et al., {\it Proceedings of the ICM, Madrid, Spain, August 22-30, 2006. Volume II: Invited lectures}.
Zurich, European Math. Soc., 2006, pp. 827--852.
\bibitem{Mi4} G. Mikhalkin. Quantum indices and refined enumeration of real plane curves.
{\it Acta Math.} {\bf 219} (2017), no. 1, 135--180.
\bibitem{NS} T. Nishinou and B. Siebert. Toric degenerations of toric
varieties and tropical curves. {\it Duke Math. J.} {\bf135} (2006), no. 1,
1--51.
\bibitem{SS} F. Schroeter and E. Shustin. Refined elliptic tropical invariants.
{\it Israel J. Math.} {\bf225} (2018),  no. 2, 817--869.
\bibitem{SpS} D. Speyer and B. Sturmfels. The tropical Grassmannian.
{\it Adv. Geom.} {\bf 4} (2004), 389--411.
\end{thebibliography}
\end{document}